# Torsion algebraic cycles and complex cobordism

Burt Totaro [*]

Atiyah and Hirzebruch gave the first counterexamples to the Hodge conjecture with integer coefficients. In particular, there is a smooth complex projective variety $X$ of dimension 7 and a torsion element of $H^4(X, \mathbf{Z})$ which is not the class of a codimension-2 algebraic cycle [**4**]. In this paper, we provide a more systematic explanation for their examples: for every smooth complex algebraic variety $X$, we show that the cycle map, from the ring of cycles modulo algebraic equivalence on $X$ to the integer cohomology of $X$, lifts canonically to a more refined topological invariant of $X$, the ring $MU^*X \otimes_{MU^*} \mathbf{Z}$, where $MU^*X$ is the complex cobordism ring of $X$. Here $MU^*X$ is a module over the graded ring $MU^* = MU^*(\text{point}) = \mathbf{Z}[x_1, x_2, \ldots]$, $x_i \in MU^{-2i}$, and we map $MU^*$ to $\mathbf{Z}$ by sending all the generators $x_i$ to 0. The ring $MU^*X \otimes_{MU^*} \mathbf{Z}$ is the same as the integer cohomology ring if the integer cohomology is torsion-free, but in general the map $MU^*X \otimes_{MU^*} \mathbf{Z} \to H^*(X, \mathbf{Z})$ need not be either injective or surjective, although the kernel and cokernel are torsion. This more refined cycle map gives a new way to prove that the Griffiths group (the kernel of the map from cycles modulo algebraic equivalence to integer cohomology) can be nonzero, without any use of Hodge theory. The resulting examples answer some questions on algebraic cycles by Colliot-Thélène and Schoen.

Our examples are all quotients of complete intersections by finite groups, as are Atiyah-Hirzebruch's examples. First, we find smooth complex projective varieties $X$ of dimension 7, definable over $\mathbf{Q}$, such that the map $CH^2X/2 \to H^4(X, \mathbf{Z}/2)$ is not injective. Here $CH^iX$ is the group of codimension $i$ algebraic cycles on $X$ modulo rational equivalence. Kollár and van Geemen [**5**], p. 135, gave the first examples of smooth complex projective varieties with $CH^2X/n \to H^4(X, \mathbf{Z}/n)$ not injective for some $n$, answering a question by Colliot-Thélène [**9**], p. 14. Over non-algebraically closed fields $k$ there are other examples of smooth projective varieties $X_k$ with $CH^2(X_k)/n \to H^4_{et}(X_k, \mathbf{Z}/n)$ not injective, due to Colliot-Thélène and Sansuc as reinterpreted by Salberger (see [**10**] and [**9**], Remark 7.6.1), Parimala and Suresh [**29**], and Bloch and Esnault [**8**]. Of these examples, only Bloch and Esnault's elements of $CH^2(X_k)/n$ are shown to remain nonzero in $CH^2(X_{\overline{k}})/n$, as happens in our example.

Also, we find codimension-3 cycles, on certain smooth complex projective varieties $X$ of dimension 15, which are torsion in the Chow group $CH^3X$, which map to 0 in $H^6(X, \mathbf{Z})$ and even in Deligne cohomology (i.e., the intermediate Jacobian), but which are not algebraically equivalent to 0. The variety $X$ and the cycles we consider can be defined over $\mathbf{Q}$. By contrast, for all $X$ over $\mathbf{C}$, the map from the torsion subgroup of $CH^iX$ to Deligne cohomology was known to be injective for $i \leq 2$ by Merkur'ev-Suslin [**22**], p. 338, and for $i = \dim X$ by Roitman [**34**], and Schoen [**35**], p. 13, asked whether the map was injective in general. Similarly, it is conjectured that codimension-2 cycles which map to 0 in Deligne cohomology are algebraically equivalent to 0, and our construction shows that this is false in codimension 3, as Nori was the first to show [**26**]. (Nori's cycles are non-torsion in the group of cycles modulo algebraic equivalence, in contrast to ours.)

We now describe the argument in more detail. In general, if $X$ is a complex algebraic scheme which can be singular or noncompact, the usual cycle class map sends the group $Z_i^{\text{alg}}X$ of $i$-dimensional algebraic cycles modulo algebraic equivalence to the Borel-Moore homology $H_{2i}^{\text{BM}}(X, \mathbf{Z})$ [**13**], chapter 19, and we lift this map to the degree $2i$ subgroup of the graded group $MU_*^{\text{BM}}X \otimes_{MU_*} \mathbf{Z}$; we informally call that subgroup $MU_{2i}^{\text{BM}}X \otimes_{MU_*} \mathbf{Z}$. Here, for a locally compact space $X$, $MU_*^{\text{BM}}X$ denotes the Borel-Moore version of the complex bordism groups of $X$. That is, for $X$ compact they are the usual complex bordism groups of $X$, generated by continuous maps of closed manifolds with a complex structure on their stable tangent bundle to $X$; and if $X = \overline{X} - S$ for a $CW$ complex $\overline{X}$ with a closed subcomplex $S$, then $MU_*^{\text{BM}}X = MU_*(\overline{X}, S)$. The construction of the new cycle map uses Hironaka's resolution of singularities together with some fundamental results on complex cobordism proved by Quillen and Wilson.

If $X$ is a locally compact space such that $H_*^{\text{BM}}(X, \mathbf{Z})$ is torsion-free, then the natural map $MU_*^{\text{BM}}X \otimes_{MU_*} \mathbf{Z} \to H_*^{\text{BM}}(X, \mathbf{Z})$ is an isomorphism; so the new cycle class map says nothing new for varieties $X$ with torsion-free homology. Also, the map $MU_*^{\text{BM}}X \otimes_{MU_*} \mathbf{Z} \to H_*^{\text{BM}}(X, \mathbf{Z})$ is always an isomorphism after tensoring with $\mathbf{Q}$. However, for spaces $X$ with torsion in their homology, the map is often not surjective, as one can compute from the fact that all odd-dimensional elements of the Steenrod algebra vanish on the image of $MU_*^{\mathrm{BM}}X$ in $H_*^{\mathrm{BM}}(X,\mathbf{Z}/p)$, for each prime number $p$. (Equivalently, the two-sided ideal in the Steenrod algebra generated by the Bockstein vanishes on the image of $MU_*^{\mathrm{BM}}X$.)

The fact that the usual cycle class map goes into the image of $MU_*^{\mathrm{BM}}X$ was essentially observed by Atiyah [**3**], footnote 1, p. 445; it follows immediately from Hironaka's resolution of singularities. As a result, all odd-dimensional elements of the Steenrod algebra vanish on the image of algebraic cycles in $H_*^{\mathrm{BM}}(X,\mathbf{Z}/p)$. Atiyah and Hirzebruch used a weaker form of this statement to produce their counterexamples to the Hodge conjecture with integer coefficients [**4**].

Once we have our more refined cycle class, it is natural to try to use it to give a new construction of smooth projective varieties with nonzero Griffiths group (the kernel of $Z^*_{\mathrm{alg}}X \to H^*(X,\mathbf{Z})$). We need to find varieties $X$ such that the map $MU^*X \otimes_{MU^*} \mathbf{Z} \to H^*(X,\mathbf{Z})$ is not injective, and then we have to hope that some of the elements of the kernel can be represented by algebraic cycles. Unfortunately, it is much harder to find topological spaces $X$ with $MU^*X \otimes_{MU^*} \mathbf{Z} \to H^*(X,\mathbf{Z})$ not injective than to find non-surjective examples. As far as I know, the only examples in the literature are those produced by Conner and Smith [**11**], [**12**], and [**37**], p. 854. Unfortunately, their examples are defined as cell complexes with explicit attaching maps, and there is no obvious way to approximate such spaces by smooth algebraic varieties.

Fortunately, there are spaces $X$ with $MU^*X \otimes_{MU^*} \mathbf{Z} \to H^*(X,\mathbf{Z})$ not injective which are more convenient for our purpose: the classifying spaces of some compact Lie groups. We will show that that the map $MU^6BG \otimes_{MU^*} \mathbf{Z} \to H^6(BG,\mathbf{Z})$ is not injective when $G$ is either the Lie group $\mathbf{Z}/2 \times SO(4)$ or the product of $\mathbf{Z}/2$ with a finite Heisenberg group $H$, a central extension

$$1 \to \mathbf{Z}/2 \to H \to (\mathbf{Z}/2)^4 \to 1.$$

There is a natural way to approximate the classifying space of a Lie group by smooth algebraic varieties, and as a result we manage to produce smooth projective varieties (the quotient of certain complete intersections by the above group $\mathbf{Z}/2 \times H$) for which our cycle map implies that the Griffiths group is not zero.

This paper was inspired by questions asked by Jack Morava. Spencer Bloch, Patrick Brosnan, Bill Fulton, Mark Green, and János Kollár had useful suggestions.

**1. A quick introduction to complex bordism.** We begin by defining a weakly complex manifold $M$ to be a smooth real manifold together with a complex vector bundle over $M$ whose underlying real vector bundle is $TM \oplus \mathbf{R}^N$ for some $N$. Thus, complex manifolds are weakly complex manifolds, but some odd-dimensional manifolds (e.g., $S^1$, since its tangent bundle is trivial) also admit weakly complex structures. We identify two complex structures on the vector bundle $TM \oplus \mathbf{R}^N$ if they are homotopic, and we also identify a complex structure on the vector bundle $TM \oplus \mathbf{R}^N$ with the obvious complex structure on $TM \oplus \mathbf{R}^N \oplus \mathbf{R}^2 = TM \oplus \mathbf{R}^{N+2}$.

The complex bordism groups $MU_iX$ of a topological space $X$, $i \geq 0$, are defined as the free abelian group on the set of continuous maps $M \to X$ where $M$ is a closed weakly complex manifold of real dimension $i$, modulo the relations

$$[M_1 \coprod M_2 \to X] = [M_1 \to X] + [M_2 \to X]$$
$$[\partial W \to X] = 0,$$

where $W$ is a compact weakly complex manifold of dimension $i+1$ with boundary together with a continuous map $W \to X$. (The boundary of $W$ inherits a weakly complex structure in a natural way.)

The notion of weakly complex manifold is rather artificial, and one might ask why we don't try to define similar invariants of a topological space $X$ using, say, complex manifolds with continuous maps to $X$. The justification for the above definition is that the groups $MU_*X$ have excellent formal properties: they form a generalized homology theory [**2**], which means that they satisfy all the usual formal properties of ordinary homology (Mayer-Vietoris, etc.) except for the dimension axiom: $MU_iX$ can be nonzero for $i > \dim X$. In fact, $MU_iX$ is always nonzero for all even $i \geq 0$, at least, because $MU_i(\mathrm{point})$ is nonzero for all even $i \geq 0$. The groups $MU_* := MU_*(\mathrm{point})$ form a ring, the product corresponding to taking products of weakly complex manifolds, and this ring was computed by Milnor and Novikov [**23**], [**27**]:

$$MU_* = \mathbf{Z}[x_1, x_2, x_3, \ldots], \qquad x_i \in MU_{2i}.$$

It happens that all the generators $x_i$ can be represented by complex manifolds. If we tensor the ring with $\mathbf{Q}$, we can take the generators to be $\mathbf{CP}^1$, $\mathbf{CP}^2$, and so on; to get the generators over $\mathbf{Z}$, we have to use certain hypersurfaces, as Milnor showed.

There is a natural map $MU_i X \to H_i(X, \mathbf{Z})$, which sends a bordism class $[M \to X]$ to the image under this map of the fundamental homology class of $M$ (since $M$ has a weakly complex structure, it has a natural orientation). This map clearly has an enormous kernel, but there is a way to define groups related to $MU_* X$ which are much closer to $H_*(X, \mathbf{Z})$. This uses that the groups $MU_* X$ form a module over the ring $MU_*$. Geometrically, the product $MU_i \otimes_\mathbf{Z} MU_j X \to MU_{i+j} X$ sends a weakly complex manifold $M^i$ and a map $M^j \to X$ to the composition $M^i \times M^j \to M^j \to X$, where the first map is the obvious projection. The point is that as long as $i > 0$, the resulting element of $MU_{i+j} X$ maps to 0 in $H_{i+j}(X, \mathbf{Z})$. So we have a natural map

$$MU_* X / (MU_{>0} \cdot MU_* X) \to H_*(X, \mathbf{Z}),$$

or, as I prefer to write it,

$$MU_* X \otimes_{MU_*} \mathbf{Z} \to H_*(X, \mathbf{Z}),$$

where the ring $MU_*$ maps to $\mathbf{Z}$ by sending all the generators $x_i$, $i \geq 1$, to 0.

If $X$ is a compact complex algebraic scheme, possibly singular, our cycle map will take values in $MU_* X \otimes_{MU_*} \mathbf{Z}$. To include schemes which may be noncompact, we define a variant of the above groups. For any locally compact topological space $X$, let $MU_i^{\mathrm{BM}} X$ be the free abelian group on the set of proper maps $M \to X$, where $M$ is a weakly complex manifold of real dimension $i$ which may be noncompact, modulo the relations

$$[M_1 \coprod M_2 \to X] = [M_1 \to X] + [M_2 \to X]$$
$$[\partial W \to X] = 0,$$

where $W$ is a weakly complex manifold of real dimension $i+1$ with boundary which may be noncompact, together with a proper map $W \to X$. These groups can be identified with bordism groups in the more usual sense for all reasonable spaces $X$: namely, the groups $MU_*^{\mathrm{BM}} X$ are the reduced bordism groups of the one-point compactification of $X$, or more generally, if $X = \overline{X} - S$ is any compactification, $MU_*^{\mathrm{BM}} X$ is isomorphic to the relative bordism group $MU_*(\overline{X}, S)$, as is defined for any generalized homology theory. The cycle map we will define for an arbitrary complex algebraic scheme $X$ takes values in $MU_*^{\mathrm{BM}} X \otimes_{MU_*} \mathbf{Z}$.

Also, as for any generalized homology theory, there is a corresponding cohomology theory, complex cobordism $MU^* X$, which is a ring for any space $X$. If $X$ is an $n$-dimensional complex manifold, there is a Poincaré duality isomorphism $MU^i X \cong MU_{2n-i}^{\mathrm{BM}} X$. So, for an $n$-dimensional complex manifold $X$, we have a geometric description of the cobordism group $MU^i X$, as bordism classes of "real codimension $i$" weakly complex manifolds $M$ (meaning that $\dim M = 2n - i$) with proper maps $M \to X$. In a sense this suffices to describe $MU^i X$ for arbitrary spaces $X$, since at least every finite cell complex is homotopy equivalent to a complex manifold (a regular neighborhood of an embedding in $\mathbf{C}^N$). Still, it may be helpful to mention one other geometric description of cobordism: for any real manifold $X$, $MU^i X$ is the group of bordism classes of codimension $i$ real manifolds $M$ with a proper map $f \colon M \to X$ and a complex structure on the "stable normal bundle" $f^* TX - TM$ [31].

If $X$ is a compact complex manifold, the above Poincaré duality isomorphism says that $MU^i X = MU_{2n-i} X$. In particular, if $X$ is a point, we have $MU^* := MU^*(\mathrm{point}) = MU_{-*}$. That is, $MU^*$ is a polynomial ring $\mathbf{Z}[x_1, x_2, \ldots]$ with $x_i \in MU^{-2i}$.

The natural way to compute $MU^* X$, for any CW complex $X$, is by the Atiyah-Hirzebruch spectral sequence $E_2 = H^*(X, MU^*) \Longrightarrow MU^* X$. This is a fourth quadrant spectral sequence because the ring $MU^*$ is in dimensions $\leq 0$; see section 6 for a picture of a similar spectral sequence. The natural map $MU^* X \to H^*(X, \mathbf{Z})$ is the "edge map" corresponding to the top row of the spectral sequence. From the spectral sequence, we can read off several of the basic properties of the complex cobordism ring, say for a finite cell complex $X$:

If $X$ has real dimension $n$, $MU^i X$ can be nonzero only for $i \leq n$. It is nonzero for all negative even $i$, at least.

The differentials are known to be torsion. It follows that $MU^* X \otimes_\mathbf{Z} \mathbf{Q}$ is a free $MU^* \otimes_\mathbf{Z} \mathbf{Q}$-module, generated by any set of elements of $MU^* X$ which map to a basis for $H^*(X, \mathbf{Q})$. In particular, the natural map $MU^* X \otimes_{MU^*} \mathbf{Z} \to H^*(X, \mathbf{Z})$ is an isomorphism $\otimes \mathbf{Q}$.

If the integer cohomology of $X$ has no torsion, then neither does the $E_2$ term of the spectral sequence. (Here the fact that the ring $MU^* = \mathbf{Z}[x_1, x_2, \ldots]$ has no torsion is crucial.) Since the differentials are always torsion, they must be 0 in this case. Thus, if $H^*(X, \mathbf{Z})$ has no torsion, then $MU^*X$ is a free $MU^*$-module, and the natural map $MU^*X \otimes_{MU^*} \mathbf{Z} \to H^*(X, \mathbf{Z})$ is an isomorphism. The corresponding statement for homology, mentioned in the introduction, follows from the homology version of the Atiyah-Hirzebruch spectral sequence.

Finally, we need to mention Brown-Peterson cohomology, a simplification of complex cobordism which is more convenient for all calculations. Namely, for each prime number $p$ there is a cohomology theory called $BP^*X$ (it is conventional not to indicate $p$ in the notation). Its coefficient ring is the polynomial ring

$$\mathbf{Z}_{(p)}[v_1, v_2, \ldots],$$

where $\mathbf{Z}_{(p)}$ is the localization of the ring $\mathbf{Z}$ at the prime $p$, and $v_1 \in BP^{-2(p-1)}$, $v_2 \in BP^{-2(p^2-1)}$, and so on. Thus the generators of the ring $BP^*$ are much more spread out than those of $MU^*$, which makes calculations easier. But $BP^*X$ carries all the topological information of $MU^*X$, because $MU^*X \otimes_{\mathbf{Z}} \mathbf{Z}_{(p)}$ splits as a direct sum of copies of $BP^*X$ in a canonical way [**30**]. In particular, $BP^*$ is a quotient ring of $MU^* \otimes_{\mathbf{Z}} \mathbf{Z}_{(p)}$ in such a way that $BP^*X = MU^*X \otimes_{MU^*} BP^*$ for all spaces $X$, and consequently

$$MU^*X \otimes_{MU^*} \mathbf{Z}_{(p)} = BP^*X \otimes_{BP^*} \mathbf{Z}_{(p)}.$$

We will use this to translate results between complex cobordism and Brown-Peterson cohomology as convenient.

**2. Quillen's theorem.** One of the fundamental facts about complex cobordism is the following classic theorem of Quillen's [**31**].

**Theorem 2.1** *Let $X$ be a finite cell complex. Then the groups $MU^*X \otimes_{MU^*} \mathbf{Z}$ are zero in negative dimensions and equal to $H^0(X, \mathbf{Z})$ in dimension 0.*

Equivalently, $MU^*X$ is generated as an $MU^*$-module by elements of nonnegative degree.

In fact, Quillen's statement can be improved a little, and we will need part of the improved statement. Namely:

**Theorem 2.2** *Let $X$ be a finite cell complex. Then the map*

$$MU^*X \otimes_{MU^*} \mathbf{Z} \to H^*(X, \mathbf{Z})$$

*is an isomorphism in dimensions $\leq 2$ and injective in dimensions $\leq 4$.*

This is best possible. In particular, the map is not surjective in dimension 3 for $X = B(\mathbf{Z}/p)^2$ or a suitable finite skeleton thereof, and it is not injective in dimension 5 for a suitable finite skeleton of $K(\mathbf{Z}, 3) \times B\mathbf{Z}/p$, as one can see by imitating the proof of Corollary 5.3 in this paper (apply it to a finite skeleton of $K(\mathbf{Z}, 3)$ in place of $BG$). We actually only need injectivity in dimensions $\leq 2$ for this paper (in the proofs of Theorem 3.1 and Lemma 4.3), except in Remark 2, section 8.

**Proof of Theorem 2.2.** This follows by the arguments Wilson used to prove Quillen's theorem [**41**], as we now explain.

The surjectivity in dimensions $\leq 2$ is trivial. In fact, an element of $H^i(X, \mathbf{Z})$ is represented by a map $X \to K(\mathbf{Z}, i)$, and for $i \leq 2$, the Eilenberg-MacLane space $K(\mathbf{Z}, i)$ has torsion-free cohomology (for $i = 0, 1, 2$, respectively, $K(\mathbf{Z}, i)$ is the space $\mathbf{Z}$, $S^1$, $\mathbf{C}P^\infty$), which implies that $MU^i K(\mathbf{Z}, i)$ maps onto the generator of $H^i(K(\mathbf{Z}, i), \mathbf{Z})$. Pulling back to $X$ proves the desired surjectivity.

It suffices to prove injectivity after tensoring with $\mathbf{Z}_{(p)}$ for each prime number $p$. As mentioned in section 1, we have

$$MU^*X \otimes_{MU^*} \mathbf{Z}_{(p)} \cong BP^*X \otimes_{BP^*} \mathbf{Z}_{(p)},$$

where $BP$ denotes Brown-Peterson cohomology at the prime $p$. So it suffices to show that the map

$$BP^*X \otimes_{BP^*} \mathbf{Z}_{(p)} \to H^*(X, \mathbf{Z}_{(p)})$$

is injective in dimensions $\leq 4$. We recall from section 1 that the coefficient ring $BP^*$ is a polynomial ring over $\mathbf{Z}_{(p)}$ with generators $v_i \in BP^{-2(p^i-1)}$, $i \geq 1$.

Following Wilson [**41**], we use the cohomology theories $BP\langle n\rangle$. These are modules over $BP$ (so $BP\langle n\rangle^*X$ is a module over $BP^*X$), with coefficients $BP\langle n\rangle^* = \mathbf{Z}_{(p)}[v_1,\ldots,v_n]$ as a $BP^*$-module (all the $v_i$'s for $i > n$ map to 0), and with maps of cohomology theories

$$BP^*X = BP\langle\infty\rangle^*X \to \cdots \to BP\langle 1\rangle^*X \to BP\langle 0\rangle^*X = H^*(X,\mathbf{Z}_{(p)}).$$

There is a long exact sequence

$$BP\langle n\rangle^{k+2(p^n-1)}X \to BP\langle n\rangle^k X \to BP\langle n-1\rangle^k X \to BP\langle n\rangle^{k+2(p^n-1)+1}X,$$

where the first map is multiplication by $v_n$ and the second map is part of the sequence of maps above. Finally, we use Wilson's main theorem [**41**], p. 118:

**Theorem 2.3** $BP^k X \to BP\langle n\rangle^k X$ *is surjective for* $k \leq 2(p^n + p^{n-1} + \cdots + 1)$.

Now we can prove Theorem 2.2. Let $x \in BP^k X$, $k \leq 4$, such that $x$ maps to 0 in $H^k(X,\mathbf{Z}_{(p)})$. We will show that $x$ is a finite sum $x = \sum_{i>0} v_i x_i$, $x_i \in BP^{k+2(p^i-1)}X$; this is equivalent to showing that $x = 0 \in BP^k X \otimes_{BP^*} \mathbf{Z}_{(p)}$.

Consider the maps

$$BP^*X \to BP\langle n\rangle^*X \to BP\langle n-1\rangle^*X.$$

If $x$ is 0 we are done. Otherwise, let $n$ be the positive integer such that $x$ maps to 0 in $BP\langle n-1\rangle^*X$ but not in $BP\langle n\rangle^*X$. Such an $n$ exists because we are assuming that $x$ maps to 0 in $BP\langle 0\rangle^*X = H^*(X,\mathbf{Z}_{(p)})$, while for $n$ large ($k$ being fixed) we have $BP^k X = BP\langle n\rangle^k X$ since $X$ is finite.

We have a commuting diagram, where the first map in each row is multiplication by $v_n$, and the second row is an exact sequence:

$$\begin{array}{ccc} BP^{k+2(p^n-1)}X & \longrightarrow & BP^k X \\ \downarrow & & \downarrow \\ BP\langle n\rangle^{k+2(p^n-1)}X & \longrightarrow BP\langle n\rangle^k X \longrightarrow & BP\langle n-1\rangle^k X \end{array}$$

Since the image $x'$ of $x$ in $BP\langle n\rangle^k X$ maps to 0 in $BP\langle n-1\rangle^k X$, $x'$ is equal to $v_n$ times an element $x'_n$ of $BP\langle n\rangle^{k+2(p^n-1)}X$. Now since $n \geq 1$ and $k \leq 4$, we have

$$k + 2(p^n - 1) \leq 4 + 2(p^n - 1) \leq 2(p^n + p^{n-1} + \cdots + 1),$$

which is exactly what we need to apply Wilson's main theorem, Theorem 2.3 above, to show that

$$BP^{k+2(p^n-1)}X \to BP\langle n\rangle^{k+2(p^n-1)}X$$

is surjective. Let $x_n$ be an element of the first group which maps to $x'_n$. Then $x - v_n x_n$ maps to 0 in $BP\langle n\rangle^k X$. Now repeat this process using $x - v_n x_n$ in place of $x$. Since $BP^k X = BP\langle N\rangle^k X$ for $N$ sufficiently large as we have said, this process stops after a finite number of steps, so that we get our finite sum $x = \sum_{i>0} v_i x_i$. QED

## 3. The new cycle map.

**Theorem 3.1** *Let $X$ be a complex algebraic scheme. We define a homomorphism from the group $Z_i^{alg}X$ of $i$-dimensional algebraic cycles on $X$ modulo algebraic equivalence to $MU_{2i}^{BM}X \otimes_{MU_*} \mathbf{Z}$ such that the composition*

$$Z_i^{alg}X \to MU_{2i}^{BM}X \otimes_{MU_*} \mathbf{Z} \to H_{2i}^{BM}(X,\mathbf{Z})$$

*is the usual cycle class map. This homomorphism is a natural transformation on the category of proper algebraic maps.*

**Proof.** The map is defined to send an irreducible $i$-dimensional subvariety $Z \subset X$ to the class of the map $[\tilde{Z} \to Z \subset X]$ in $MU_{2i}^{\mathrm{BM}} X \otimes_{MU_*} \mathbf{Z}$, where $\tilde{Z} \to Z$ is any resolution of singularities of $Z$, that is, a proper birational map with $\tilde{Z}$ smooth. Such resolutions exist, by Hironaka [16]. The first step is to show that the various elements of $MU_{2i}^{\mathrm{BM}} X$ that can arise from different resolutions $\tilde{Z}$ of a fixed variety $Z$ are all equal in $MU_{2i}^{\mathrm{BM}} X \otimes_{MU_*} \mathbf{Z}$.

Let $\tilde{Z}_1$ and $\tilde{Z}_2$ be any two resolutions of $Z$. By Hironaka [16], it is possible to blow up $\tilde{Z}_1$ repeatedly along smooth subvarieties to get a variety $\tilde{Z}_1'$ which maps to $\tilde{Z}_2$, giving a commutative diagram.

$$\begin{array}{ccc} \tilde{Z}_1' & \longrightarrow & \tilde{Z}_1 \\ \downarrow & & \downarrow \\ \tilde{Z}_2 & \longrightarrow & Z \end{array}$$

By Quillen's theorem (Theorem 2.1 above), for any finite complex $X$, the group $MU^i X \otimes_{MU^*} \mathbf{Z}$ is 0 for $i < 0$ and equals $H^0(X, Z)$ for $i = 0$. By Poincaré duality for complex cobordism, it follows that, for any smooth complex $n$-manifold, the group $MU_i^{\mathrm{BM}} X \otimes_{MU_*} \mathbf{Z}$ is 0 for $i > 2n$ and equals $H_{2n}^{\mathrm{BM}}(X, \mathbf{Z})$ for $i = 2n$. In particular, if $X \to Y$ is a proper birational morphism of smooth $n$-dimensional complex varieties, we have

$$[X \to Y] = [Y \to Y] \in MU_{2n}^{\mathrm{BM}}(Y, \mathbf{Z}) \otimes_{MU_*} \mathbf{Z},$$

because this is true in $H_{2n}^{\mathrm{BM}}(Y, \mathbf{Z}) = \mathbf{Z}$. Thus, in the situation of the previous paragraph, we have

$$[\tilde{Z}_1' \to \tilde{Z}_1] = [\tilde{Z}_1 = \tilde{Z}_1] \in MU_{2n}^{\mathrm{BM}} \tilde{Z}_1 \otimes_{MU_*} \mathbf{Z}$$

and

$$[\tilde{Z}_1' \to \tilde{Z}_2] = [\tilde{Z}_2 = \tilde{Z}_2] \in MU_{2n}^{\mathrm{BM}} \tilde{Z}_2 \otimes_{MU_*} \mathbf{Z}.$$

It follows that

$$[\tilde{Z}_1 \to X] = [\tilde{Z}_1' \to X] = [\tilde{Z}_2 \to X]$$

in $MU_{2n}^{\mathrm{BM}} X \otimes_{MU_*} \mathbf{Z}$. Thus, any two resolutions of a subvariety $Z \subset X$ define the same element of $MU_{2n}^{\mathrm{BM}} X \otimes_{MU_*} \mathbf{Z}$, which we are now justified in calling the class $[Z]$ of $Z$ in $MU_{2n}^{\mathrm{BM}} X \otimes_{MU_*} \mathbf{Z}$.

Thus we have a natural map $Z_n X \to MU_{2n}^{\mathrm{BM}} X \otimes_{MU_*} \mathbf{Z}$ for any complex algebraic scheme $X$, where $Z_n X$ is the group of algebraic $n$-cycles on $X$, that is, the free abelian group on the set of closed $n$-dimensional irreducible subvarieties of $X$. For a cycle $\alpha$, we write $[\alpha]$ for its class in $MU_*^{\mathrm{BM}} X \otimes_{MU_*} \mathbf{Z}$. Our next step will be to check that this map is a natural transformation on the category of proper algebraic maps.

We recall the definition of the map $f_* : Z_n X \to Z_n Y$ associated to a proper algebraic map $f : X \to Y$ [13]. For a closed subvariety $Z \subset X$, $f(Z)$ is a closed subvariety of $Y$, and we define

$$f_*(Z) = \begin{cases} \deg\,(f : Z \to f(Z)) f(Z) & \text{if} \quad \dim f(Z) = \dim Z \\ 0 & \text{if} \quad \dim f(Z) < \dim Z. \end{cases}$$

So let $X \to Y$ be a proper algebraic map, and $Z$ an $n$-dimensional subvariety of $X$. To show that the map $Z_n X \to MU_{2n}^{\mathrm{BM}} X \otimes_{MU_*} \mathbf{Z}$ is a natural transformation means to show that the class in $MU_{2n}^{\mathrm{BM}} Y \otimes_{MU_*} \mathbf{Z}$ of the cycle $f_*(Z)$ is equal to the image of the class $[Z] \in MU_{2n}^{\mathrm{BM}} X \otimes_{MU_*} \mathbf{Z}$ under the natural map

$$f_* : MU_*^{\mathrm{BM}} X \to MU_*^{\mathrm{BM}} Y.$$

There are two cases, depending on whether the dimension of $f(Z)$ is $n$ or less than $n$.

If $f(Z)$ has dimension $n$, let $\tilde{Z}_2$ be a resolution of singularities of $f(Z)$, and let $\tilde{Z}_1$ be a resolution of singularities of $Z$ such that the rational map from $Z$ to $\tilde{Z}_2$ becomes well-defined on $\tilde{Z}_1$.

$$\begin{array}{ccccc} \tilde{Z}_1 & \longrightarrow & Z & \longrightarrow & X \\ \downarrow & & \downarrow & & \downarrow \\ \tilde{Z}_2 & \longrightarrow & f(Z) & \longrightarrow & Y \end{array}$$

Clearly the map $\tilde{Z}_1 \to \tilde{Z}_2$ has the same degree $d$ as the map $f : Z \to f(Z)$. By Quillen's theorem, Theorem 2.1 above, we have
$$[\tilde{Z}_1 \to \tilde{Z}_2] = d[\tilde{Z}_2 = \tilde{Z}_2] \in MU_{2n}^{BM}\tilde{Z}_2 \otimes_{MU_*} \mathbf{Z},$$
since $MU_{2n}^{BM}\tilde{Z}_2 \otimes_{MU_*} \mathbf{Z} \cong H_{2n}^{BM}(\tilde{Z}_2, \mathbf{Z}) = \mathbf{Z}$. Now the image under $f_* : MU_*^{BM}X \otimes_{MU_*} \mathbf{Z} \to MU_*^{BM}Y \otimes_{MU_*} \mathbf{Z}$ of the class $[Z]$ is by definition equal to the class of the map $\tilde{Z}_1 \to Y$, and since this map factors through $\tilde{Z}_2$, the above equality means that this class is equal to $d$ times $[\tilde{Z}_2 \to Y]$, that is, to $d$ times $[Z_2]$. This proves functoriality of our map in this case.

The argument is similar if $f(Z)$ has dimension less than $n$. In this case we use Quillen's theorem to prove that a proper holomorphic map $X \to Y$ between complex manifolds with $\dim X > \dim Y$ has $[X \to Y] = 0 \in MU_*^{BM}Y \otimes_{MU_*} \mathbf{Z}$, because this group is 0. As in the previous case, we apply this result to a resolution $\tilde{Z}_2$ of $f(Z)$ and a resolution $\tilde{Z}_1$ of $Z$ which maps to $\tilde{Z}_2$, and we find that $f_*[Z_1] = 0 \in MU_*^{BM}Y \otimes_{MU_*} \mathbf{Z}$.

Thus we have defined a natural transformation $Z_*X \to MU_*^{BM}X \otimes_{MU_*} \mathbf{Z}$ on the category of complex algebraic schemes and proper algebraic maps. To finish the proof of the theorem, we have to show that this map is well defined on cycles modulo algebraic equivalence. That is, we have to show that for every smooth compact connected curve $C$ and every $n+1$-dimensional subvariety $W \subset X \times C$ with the second projection $f : W \to C$ not constant, we have $[(p_1)_* f^*(a)] = [(p_1)_* f^*(b)] \in MU_{2n}^{BM}X \otimes_{MU_*} \mathbf{Z}$ for every pair of points $a, b \in \mathbf{C}$, where $p_1 : W \to X$ is the first projection. (The cycles $f^*(a)$ and $f^*(b)$ in $W$ are defined in [**13**], Chapter 2, as the Weil divisors associated to the obvious Cartier divisors.) In view of the naturality we have proved, it suffices to prove that $[f^*(a)] = [f^*(b)]$ in $MU_{2n}^{BM}W \otimes_{MU_*} \mathbf{Z}$.

Let $\pi : \tilde{W} \to W$ be a resolution of singularities of $W$. The pushforwards of the cycles $(f\pi)^*(a)$, $(f\pi)^*(b)$ on $\tilde{W}$ to $W$ are the cycles $f^*(a)$, $f^*(b)$ on $W$, by Fulton [**13**], p. 34, proof of (c). So it suffices to prove that $[(f\pi)^*(a)] = [(f\pi)^*(b)] \in MU_{2n}^{BM}\tilde{W} \otimes_{MU_*} \mathbf{Z}$. But we know that algebraically equivalent cycles are homologous, so that these two cycles are equal in $H_{2n}^{BM}(\tilde{W}, \mathbf{Z})$; and by the extension of Quillen's theorem given in Theorem 2.2, we have $MU^2\tilde{W} \otimes_{MU^*} \mathbf{Z} \cong H^2(\tilde{W}, \mathbf{Z})$, that is, by Poincaré duality on the smooth $(n+1)$-dimensional variety $\tilde{W}$, $MU_{2n}^{BM}\tilde{W} \otimes_{MU_*} \mathbf{Z} = H_{2n}^{BM}(\tilde{W}, \mathbf{Z})$. Thus $[(f\pi)^*(a)] = [(f\pi)^*(b)] \in MU_{2n}^{BM}\tilde{W} \otimes_{MU_*} \mathbf{Z}$, and so the cycle map is well-defined on algebraic equivalence classes. QED

**Remark.** The cycle class map is in fact well-defined on a slightly weaker equivalence relation than algebraic equivalence, as explained in Remark 2, section 8.

**4. Products.** If $X$ is a smooth complex algebraic variety of dimension $n$, then cycles modulo algebraic equivalence, graded by $Z_{alg}^i X = Z_{n-i}^{alg} X$, form a ring, as do the groups $MU^i X \otimes_{MU^*} \mathbf{Z} = MU_{2n-i}^{BM}X \otimes_{MU_*} \mathbf{Z}$.

**Theorem 4.1** *If $X$ is a smooth variety, then the cycle map $Z_{alg}^i X \to MU^{2i}X \otimes_{MU^*} \mathbf{Z}$ is a ring homomorphism.*

The proof follows the outline of Fulton's proof that the usual cycle map $Z_{alg}^* X \to H^*(X, \mathbf{Z})$ is a ring homomorphism [**13**], chapter 19.

**Proof.** It is equivalent to check that the map from cycles modulo rational equivalence to $MU^*X \otimes_{MU^*} \mathbf{Z}$ is a ring homomorphism. (We do things this way because most of Fulton's book is written in terms of cycles modulo rational equivalence; the same arguments would apply to cycles modulo algebraic equivalence.)

We recall the construction of the intersection product on cycles modulo rational equivalence given by Fulton and MacPherson [**13**]. Given cycles $\alpha$ and $\beta$ on a smooth variety $X$, there is a product cycle $\alpha \times \beta$ on $X \times X$, and the product $\alpha\beta \in CH_*X$ is defined as the pullback of $\alpha \times \beta$ to the diagonal by a map $CH_i(X \times X) \to CH_{i-n}X$. This pullback map is defined, more generally, for any regular embedding: if $X$ is any local complete intersection subscheme of codimension $d$ in a scheme $Y$, then there is a pullback map $CH_iY \to CH_{i-d}X$. For the fundamental example of a regular embedding, the inclusion of the zero-section of a vector bundle into the total space of the vector bundle, $X \hookrightarrow E$, the pullback map $CH_iE \to CH_{i-d}X$ ($d$ = rank $E$) is defined to be the inverse of the natural map $CH_{i-d}X \to CH_iE$, sending a subvariety $Z \subset X$ to $E|_Z \subset E$, which one proves to be an isomorphism. For an arbitrary regular embedding $X \to Y$ of codimension $d$, the pullback map sends a subvariety $V \subset Y$ to the pullback under the zero-section inclusion $X \to N_{X/Y}$ of the normal cone $C$ to $V \cap X$ in $V$. Here the normal cone $C$ (defined as Spec $(\oplus_{n \geq 0} \mathcal{I}^n/\mathcal{I}^{n+1})$, where $\mathcal{I}$ is the ideal sheaf defining $V \cap X$ in $V$) is a subscheme of the normal bundle $N_{X/Y}$ of the same dimension, $i$, as $V$, so $C$ gives an $i$-dimensional cycle in $N_{X/Y}$, which pulls back to an element of $CH_{i-d}X$ by the map we have already defined.

The product on $MU_*^{\mathrm{BM}}X$, for a complex $n$-manifold $X$, can be defined similarly. There is an external product $MU_*^{\mathrm{BM}}X \otimes_{MU^*} MU_*^{\mathrm{BM}}X \to MU_*^{\mathrm{BM}}(X \times X)$, and the internal product is defined by composing that with a pullback map $MU_i^{\mathrm{BM}}(X \times X) \to MU_{i-2n}^{\mathrm{BM}}X$. The pullback map is defined more generally: for any codimension-$d$ complex submanifold $X$ of a complex manifold $Y$, there is a pullback map $MU_i^{\mathrm{BM}}Y \to MU_{i-2d}^{\mathrm{BM}}X$. It can be defined as cap product with an "orientation class" $u_{XY} \in MU^{2d}(Y, Y-X)$. To define $u_{XY}$, identify $MU^{2d}(Y, Y-X)$ with $MU^{2d}(N_{X/Y}, N_{X/Y}-X)$ by excision (where $X$ is included in the normal bundle $N_{X/Y}$ as the zero-section); then $u_{XY}$ is the Thom class of the complex vector bundle $N_{X/Y}$. (This pullback map is easier to define than the one on Chow groups, because a tubular neighborhood of $X \subset Y$ is diffeomorphic to the normal bundle $N_{X/Y}$, whereas in algebraic geometry there is typically no neighborhood of $X \subset Y$ which is algebraically or even analytically isomorphic to the normal bundle of $X$.)

Comparing these constructions, we see that the theorem would follow from the commutativity of the diagram of pullback maps,

$$\begin{array}{ccc} CH_i Y & \longrightarrow & CH_{i-d} X \\ \downarrow & & \downarrow \\ MU_{2i}^{\mathrm{BM}} Y \otimes_{MU_*} \mathbf{Z} & \longrightarrow & MU_{2(i-d)}^{\mathrm{BM}} X \otimes_{MU_*} \mathbf{Z}, \end{array}$$

for any codimension-$d$ smooth algebraic subvariety $X$ of a smooth algebraic variety $Y$. It seems natural to prove this in somewhat greater generality. Namely, for any codimension-$d$ regular embedding $X \to Y$ of complex algebraic schemes, Baum, Fulton, and MacPherson define an orientation class $u_{XY} \in MU^{2d}(Y, Y - X)$ which agrees with that defined above when $X$ and $Y$ are smooth [**6**], p. 137. (In fact, they define such a class in any complex-oriented cohomology theory, that is, any generalized cohomology theory with Thom classes for complex vector bundles. Complex cobordism is complex-oriented; in fact, it is the universal such theory, so that the Baum-Fulton-MacPherson orientation class in complex cobordism maps to the corresponding orientation class in any other complex-oriented cohomology theory.) Briefly, they extend the normal bundle $N_{X/Y}$ to a topological complex vector bundle $Q$ on a tubular neighborhood $N$ of $X$, and they construct a continuous section of $Q$ which vanishes on $X$ "as a scheme"; then this section gives a map $(N, N - X) \to (Q, Q - N)$, and $u_{XY}$ is the pullback of the Thom class in $MU^{2d}(Q, Q - N)$ to $MU^{2d}(N, N - X) = MU^{2d}(Y, Y - X)$. Cap product with this class $u_{XY}$ defines a pullback map $MU_i^{\mathrm{BM}} Y \to MU_{i-2d}^{\mathrm{BM}} X$ of $MU_*$-modules.

The theorem follows from the following statement about that pullback map.

**Lemma 4.2** *For any codimension-$d$ regular embedding of complex algebraic schemes $X \to Y$, the following diagram of pullback maps commutes.*

$$\begin{array}{ccc} CH_i Y & \longrightarrow & CH_{i-d} X \\ \downarrow & & \downarrow \\ MU_{2i}^{BM} Y \otimes_{MU_*} \mathbf{Z} & \longrightarrow & MU_{2(i-d)}^{BM} X \otimes_{MU_*} \mathbf{Z} \end{array}$$

**Proof.** This is easy to check for the fundamental example of regular embeddings, the inclusion $X \to E$ of the zero-section of a vector bundle over a scheme $X$. In this case, the pullback maps $CH_i E \to CH_{i-d} X$ and $MU_i^{\mathrm{BM}} E \to MU_{i-2d}^{\mathrm{BM}} X$ are both isomorphisms, so it suffices to prove that the inverse maps commute:

$$\begin{array}{ccc} CH_{i-d} X & \longrightarrow & CH_i E \\ \downarrow & & \downarrow \\ MU_{2(i-d)}^{\mathrm{BM}} X \otimes_{MU_*} \mathbf{Z} & \longrightarrow & MU_{2i}^{\mathrm{BM}} E \otimes_{MU_*} \mathbf{Z}. \end{array}$$

The map on the top row sends a subvariety $Z \subset X$ to the subvariety $E|_Z \subset E$, and the map on the bottom row sends a proper map $\tilde{Z} \to X$, for a weakly complex manifold $\tilde{Z}$, to the obvious proper map $E|_{\tilde{Z}} \to E$.

The cycle map sends $Z \subset X$ to the map $\tilde{Z} \to X$ for any resolution $\tilde{Z}$ of $X$; since $E|_{\tilde{Z}}$ is a resolution of $E|_Z$, the diagram commutes.

Also, we need to check the lemma for regular embeddings of codimension one, at least in the following situation.

**Lemma 4.3** *Let $X$ be an $n$-dimensional variety, $T$ a smooth curve, $f : X \to T$ a surjective map, $t \in T$. Then*
$$f^* u_{t,T} \cap [X] = [f^*(t)] \in MU^{BM}_{2n-2}(f^{-1}(t)) \otimes_{MU_*} \mathbf{Z},$$
*where $u_{t,T} \in MU^2(T, T - t)$ is the orientation class and $f^*(t)$ is the Weil divisor associated to the obvious Cartier divisor.*

**Proof.** This is similar to the proof that the cycle map is well-defined on algebraic equivalence (Theorem 3.1). Namely, let $\pi : \tilde{X} \to X$ be a resolution. By Fulton [**13**], p. 34, proof of (c), we have $\pi_*((f\pi)^*(t)) = f^*(t) \in Z_{n-1} f^{-1}(t)$. Also, in $MU^{BM}_{2n-2}(f^{-1}(t)) \otimes_{MU_*} \mathbf{Z}$, we have

$$\begin{aligned}
\pi_*((f\pi)^* u_{t,T} \cap [\tilde{X}]) &= \pi_*(\pi^* f^* u_{t,T} \cap [\tilde{X}]) \\
&= f^* u_{t,T} \cap \pi_*[\tilde{X}] \\
&= f^* u_{t,T} \cap [X].
\end{aligned}$$

As a result, it suffices to show that
$$f^* u_{t,T} \cap [X] = [f^*(t)] \in MU^{BM}_{2n-2}(f^{-1}(t)) \otimes_{MU_*} \mathbf{Z}$$
for $X$ smooth of dimension $n$.

For $X$ smooth, we can identify $MU^{BM}_{2n-2}(f^{-1}(t)) \otimes_{MU_*} \mathbf{Z}$ with $MU^2(X, X - f^{-1}(t)) \otimes_{MU^*} \mathbf{Z}$. This injects into $H^2(X, X - f^{-1}(t), \mathbf{Z})$ by the extension of Quillen's theorem given in Theorem 2.2. (The theorem is stated in terms of the cohomology of a space, rather than a pair of spaces, but for any generalized cohomology theory $h^*$ we can identify $h^i(X, X - S)$ with reduced $h^i$ of a pointed space (the mapping cylinder modulo $X - S$), and the theorem clearly applies to the reduced cohomology of a pointed space.) So it suffices to prove the above equality in $H^2(X, X - f^{-1}(t), \mathbf{Z}) = H^{BM}_{2n-2}(f^{-1}(t), \mathbf{Z})$, as is done in [**13**], p. 373. QED

Now we can prove Lemma 4.2 for a general regular embedding $X \to Y$. As mentioned at the beginning of the proof of Theorem 4.1, the pullback map $CH_i Y \to CH_{i-d} X$ sends a subvariety $V \subset Y$ to the pullback of $C_{V \cap X} V \subset N_{X/Y}$ to $X$. Since we have checked the lemma for the inclusion of $X$ into the vector bundle $N_{X/Y}$, the lemma in general reduces to the statement that the pullback of $[V] \in MU^{BM}_{2i} Y \otimes_{MU_*} \mathbf{Z}$ to $X$ is equal to the pullback of the class of the normal cone $[C_{V \cap X} V] \in MU^{BM}_{2i} N_{X/Y} \otimes_{MU_*} \mathbf{Z}$ to $X$.

To prove this, we use that every embedding of one scheme in another has a natural "deformation to the normal cone" [**13**], chapter 5. That is, for a closed subscheme $X \subset Y$, there is a scheme $M_X Y$ and a $\mathbf{P}^1$-family of embeddings of $X$, $X \times \mathbf{P}^1 \subset M_X Y$, with a commutative diagram

$$\begin{array}{ccc} X \times \mathbf{P}^1 & \longrightarrow & M_X Y \\ & \searrow & \downarrow \rho \\ & & \mathbf{P}^1, \end{array}$$

such that for $t \in \mathbf{P}^1 - \{\infty\} = A^1$, $\rho^{-1}(t) \cong Y$ and the embedding $X \subset \rho^{-1}(t)$ is the given embedding $X \subset Y$, and over $\infty$, the embedding is the embedding of $X$ in the normal cone $C_X Y$ of $X$ in $Y$. The map $M_X Y \to \mathbf{P}^1$ is flat, and so the inverse image of each point of $\mathbf{P}^1$ is a Cartier divisor. Explicitly, $M_X Y$ is the blow-up of $Y \times \mathbf{P}^1$ along $X \times \infty$, with the proper transform of $Y \times \infty$ omitted. It is worth mentioning that the normal cone is a vector bundle over $X$ (what we have been calling the normal bundle $N_{X/Y}$) if and only if the inclusion $X \to Y$ is a regular embedding.

To prove the lemma, let $X \to Y$ be a regular embedding and $V \subset Y$ a subvariety. Then $M' = M_{V \cap X}V$ is a subvariety of the scheme $M = M_X Y$.

$$
\begin{array}{ccccc}
V \cap X & \to & C & \to & \{\infty\} \\
\downarrow & & \downarrow & & \downarrow \\
(V \cap X) \times \mathbf{P}^1 & \to & M' & \to & \mathbf{P}^1 \\
\uparrow & & \uparrow & & \uparrow \\
V \cap X & \to & V & \to & \{0\}
\end{array}
\quad \subset \quad
\begin{array}{ccccc}
X & \to & N_{X/Y} & \to & \{\infty\} \\
\downarrow & & \downarrow & & \downarrow \\
X \times \mathbf{P}^1 & \xrightarrow{F} & M & \to & \mathbf{P}^1 \\
\uparrow & & \uparrow & & \uparrow \\
X & \to & Y & \to & \{0\}
\end{array}
$$

The variety $V$ and the normal cone $C = C_{V \cap X}V$ are the fibers over 0 and $\infty$, respectively, of the map $M' \to \mathbf{P}^1$, so Lemma 4.3 gives that the fundamental class in $MU_{2i}^{\mathrm{BM}}V \otimes_{MU_*} \mathbf{Z}$ of the subvariety $V \subset Y$ is the Baum-Fulton-MacPherson pullback of the fundamental class of $M'$, by the codimension-one regular embedding $V \subset M'$. It follows that that the class $[V] \in MU_{2i}^{\mathrm{BM}}Y \otimes_{MU_*} \mathbf{Z}$ is the pullback of the class $[M'] \in MU_{2i}^{\mathrm{BM}}M \otimes_{MU_*} \mathbf{Z}$ by the codimension-one regular embedding $Y \subset M$. Likewise, the class in $MU_{2i}^{\mathrm{BM}}N_{X/Y} \otimes_{MU_*} \mathbf{Z}$ of the normal cone $C = C_{V \cap X}V$, viewed as a scheme and thus as a cycle with multiplicities, is the pullback of the class of $[M']$ on $M$ by the codimension-one regular embedding $N_{X/Y} \subset M$.

As a result, the pullback of the class $[V]$ on $Y$ to $MU_{2(i-d)}^{\mathrm{BM}}X \otimes_{MU_*} \mathbf{Z}$ is the pullback of $[M']$ on $M$ first by the codimension-one regular embedding $Y \subset M$, then by the codimension-$d$ regular embedding $X \times 0 \subset Y$. Likewise, the pullback of the class $[C]$ on $N_{X/Y}$ to $X$ is the pullback of the same class $[M']$ first by the codimension-one regular embedding $N_{X/Y} \subset M$, then by the codimension-$d$ regular embedding $X \subset N_{X/Y}$. By the above commutative diagrams, using the naturality of Baum-Fulton-MacPherson's pullback maps on bordism groups, we find that these two elements of $MU_{2(i-d)}^{\mathrm{BM}}X \otimes_{MU_*} \mathbf{Z}$ are the pullbacks of the same class $[M']$ on $M$ by the codimension-$d$ regular embedding $F: X \times \mathbf{P}^1 \to M$, followed by the pullback to $X \times 0$ or $X \times \infty$ respectively. Since these last two pullbacks are equal on $MU_*^{\mathrm{BM}}(X \times \mathbf{P}^1)$, the two elements of $MU_{2(i-d)}^{\mathrm{BM}}X \otimes_{MU_*} \mathbf{Z}$ are equal. QED

**5. Non-injectivity of the map $MU^*X \otimes_{MU^*} \mathbf{Z} \to H^*(X, \mathbf{Z})$.** In this section we construct some topological spaces $X$ for which the map $MU^*X \otimes_{MU^*} \mathbf{Z} \to H^*(X, \mathbf{Z})$ is not injective, and such that there is a natural way to approximate the homotopy type of $X$ by smooth algebraic varieties. Namely, $X$ will be the classifying space $BG$ of a compact Lie group $G$ (we will explain the relation to algebraic geometry in section 7).

We need to take several precautions when talking about the complex cobordism of an infinite CW complex such as $BG$. The point is that $MU^*(\text{point})$ is nonzero in all even dimensions $\leq 0$, so that $MU^iX$ is affected by all the cells in $X$ of dimension $\geq i$. One phenomenon here is described by Milnor's exact sequence, which holds for all generalized cohomology theories $h^*$ and all infinite CW complexes $X$ [**24**]:

$$0 \to \varprojlim{}^1_n h^{i-1}(X_n) \to h^i X \to \varprojlim{}_n h^i(X_n) \to 0.$$

Here $X_n$ denotes the $n$-skeleton of $X$. The Atiyah-Hirzebruch spectral sequence $H^*(X, h^*) \Longrightarrow h^*(X)$ actually converges to $\varprojlim h^*(X_n)$, not to $h^*X$.

Fortunately, the only infinite complexes we will need to consider are classifying spaces of compact Lie groups $BG$, in which case the $\varprojlim{}^1$ term for $MU$-theory and the other cohomology theories we consider is 0, by Landweber [**21**]. For such spaces, Landweber also proves a strong Mittag-Leffler statement about the inverse limit in $MU^*X = \varprojlim MU^*X_n$: namely, for each $n$ there is an $m \geq n$ such that (in all dimensions at once) we have im $(MU^*X \to MU^*X_n) = $ im $(MU^*X_m \to MU^*X_n)$. To give a little context for these statements: for the space $X = K(\mathbf{Z}, 3)$, which is outside the class we consider, the $\varprojlim{}^1$ group is nonzero, and the Mittag-Leffler statement fails, because im $(MU^3K(\mathbf{Z},3)_n \to H^3(K(\mathbf{Z},3)_n, \mathbf{Z}) = \mathbf{Z})$ is a subgroup of finite index which decreases to 0 as $n$ goes to infinity. See [**33**] for some clarification of this phenomenon.

Since each group $MU^iBG$ is an inverse limit, we have to view it as a topological abelian group. In particular, tensor products involving $MU^*BG$ will always mean *completed* tensor products. For example,

$$MU^iBG \otimes_{MU^*} \mathbf{Z} := \varprojlim \quad [(\text{im } MU^*BG \to MU^*(BG)_n) \otimes_{MU^*} \mathbf{Z}].$$

Likewise, following Kono and Yagita [19], we define

$$MU^*BG \otimes_{MU^*} MU^*BH := \varprojlim \quad (\text{im } MU^*BG \to MU^*(BG)_n) \otimes_{MU^*} (\text{im } MU^*BH \to MU^*(BH)_n).$$

Now we can state the main result of this section.

**Theorem 5.1** *Let $G$ be either $SO(4)$ or the central extension $1 \to \mathbf{Z}/2 \to G \to (\mathbf{Z}/2)^4 \to 1$ contained in $SO(4)$. Then the maps $MU^4BG \otimes_{MU^*} \mathbf{Z}/2 \to H^4(BG, \mathbf{Z}/2)$ and $MU^*(BG \times B\mathbf{Z}/2) \otimes_{MU^*} \mathbf{Z} \to H^6(BG \times B\mathbf{Z}/2, \mathbf{Z})$ are not injective.*

We begin by explaining why product groups are convenient for this question. In the simplest examples, $MU^*BG \otimes_{MU^*} \mathbf{Z} \to H^*(BG, \mathbf{Z})$ tends to be injective but not surjective. It happens, however, that even if this map is injective for two groups $G_1$ and $G_2$, it need not be injective for $G_1 \times G_2$. Specifically, if the map $MU^*BG \otimes_{MU^*} \mathbf{Z} \to H^*(BG, \mathbf{Z})$ is injective but not split injective, as a map of abelian groups, then the map $MU^*(BG \times B\mathbf{Z}/p) \otimes_{MU^*} \mathbf{Z} \to H^*(BG \times B\mathbf{Z}/p, \mathbf{Z})$ is not injective for some prime $p$. The following lemma expresses this idea more precisely.

**Lemma 5.2** *Let $G$ be a compact Lie group. If the map $MU^kBG \otimes_{MU^*} \mathbf{Z}/p \to H^k(BG, \mathbf{Z}/p)$ is not injective, then $MU^{k+2}(BG \times B\mathbf{Z}/p) \otimes_{MU^*} \mathbf{Z} \to H^{k+2}(BG \times B\mathbf{Z}/p, \mathbf{Z})$ is not injective.*

**Proof.** By Landweber [20], Theorem 2$'$, $MU^*B\mathbf{Z}/p$ is flat in the sense that we have an isomorphism of topological abelian groups

$$MU^*(X \times B\mathbf{Z}/p) = MU^*X \otimes_{MU^*} MU^*B\mathbf{Z}/p$$

for all finite complexes $X$. (Kono and Yagita [19] conjecture the same statement for all compact Lie groups in place of $\mathbf{Z}/p$.) As a consequence, we have an isomorphism of topological abelian groups

$$MU^*(BG \times B\mathbf{Z}/p) = MU^*BG \otimes_{MU^*} MU^*B\mathbf{Z}/p$$

for all compact Lie groups $G$. The right side is a completed tensor product, as explained above.

It follows that

$$MU^*(BG \times B\mathbf{Z}/p) \otimes_{MU^*} \mathbf{Z} = (MU^*BG \otimes_{MU^*} \mathbf{Z}) \otimes_{\mathbf{Z}} (MU^*B\mathbf{Z}/p \otimes_{MU^*} \mathbf{Z}),$$

where the right side is again a completed tensor product. Here $MU^*B\mathbf{Z}/p \otimes_{MU^*} \mathbf{Z}$ maps isomorphically to $H^*(B\mathbf{Z}/p, \mathbf{Z}) = \mathbf{Z}[[c_1]]/(pc_1)$, where $c_1$ is in dimension 2. So

$$MU^*(BG \times B\mathbf{Z}/p) \otimes_{MU^*} \mathbf{Z} = MU^*BG \otimes_{MU^*} \mathbf{Z} \oplus \prod_{i \geq 1} c_1^i \, MU^*BG \otimes_{MU^*} \mathbf{Z}/p,$$

and by our assumption in this lemma, the group $c_1 \cdot MU^kBG \otimes_{MU^*} \mathbf{Z}/p$ maps non-injectively to the corresponding group in

$$H^*(BG, \mathbf{Z}) \oplus \prod_{i \geq 1} c_1^i \cdot H^*(BG, \mathbf{Z}/p)$$
$$= H^*(BG, \mathbf{Z}) \otimes_{\mathbf{Z}} H^*(B\mathbf{Z}/p, \mathbf{Z})$$
$$\subset H^*(BG \times B\mathbf{Z}/p, \mathbf{Z}).$$

Thus $MU^{k+2}(BG \times B\mathbf{Z}/p) \otimes_{MU^*} \mathbf{Z} \to H^{k+2}(BG \times B\mathbf{Z}/p, \mathbf{Z})$ is not injective.  QED

**Corollary 5.3** *Let $G$ be a compact Lie group. If the image of the map $MU^k BG \to H^k(BG, \mathbf{Z})$ contains $p$ times an element $x \in H^k(BG, \mathbf{Z})$, for some prime number $p$, but does not contain $x$ itself or $x$ plus any element killed by $p$ in $H^k(BG, \mathbf{Z})$, then $MU^k BG \otimes_{MU^*} \mathbf{Z}/p \to H^*(BG, \mathbf{Z}/p)$ and $MU^{k+2}(BG \times B\mathbf{Z}/p) \otimes_{MU^*} \mathbf{Z} \to H^{k+2}(BG \times B\mathbf{Z}/p, \mathbf{Z})$ are not injective.*

**Proof.** The hypothesis implies that an element of $MU^k BG$ which maps to $px \in H^k(BG, \mathbf{Z})$ is nonzero in $MU^k(BG) \otimes_{MU^*} \mathbf{Z}/p$, and it clearly maps to 0 in $H^k(BG, \mathbf{Z}/p)$. Thus Lemma 5.2 applies. QED

Thus, to prove Theorem 5.1, it suffices to prove that the image of $MU^4 BG \to H^4(BG, \mathbf{Z})$ contains 2 times some element $x \in H^4(BG, \mathbf{Z})$ but not $x$ plus any 2-torsion element, when $G$ is $SO(4)$ or the group of order 32 mentioned in the theorem. For $G = SO(4)$, Kono and Yagita [19] computed $MU^* BSO(4)$, and we can read this off from their calculation. For clarity, here is a direct proof. The point is that there is a class $\chi \in H^4(BSO(4), \mathbf{Z})$, the Euler class, such that $2\chi$ is in the image of $MU^4 BSO(4)$ but $\chi$ is not. See Milnor-Stasheff [25] for the cohomology of $BSO(n)$ and in particular the definition of the Euler class. In fact, in $H^4(BSO(4), \mathbf{Z}) = \mathbf{Z} \oplus \mathbf{Z}$, we can compute that

$$2\chi = c_2 A - c_2 B,$$

where $B: SO(4) \to GL(4, \mathbf{C})$ is the obvious representation and $A$ is the representation given by thinking of $SO(4)$ as a double cover of $SO(3) \times SO(3)$, and then projecting to the first $SO(3)$:

$$A: SO(4) \to SO(3) \times SO(3) \to SO(3) \to GL(3, \mathbf{C}).$$

(To check that $2\chi = c_2 A - c_2 B$, use that $SO(4)$ is doubly covered by $SU(2) \times SU(2)$, and $H^4(BSO(4), \mathbf{Z})$ injects into $H^4(BSU(2) \times BSU(2), \mathbf{Z}) = \mathbf{Z} \oplus \mathbf{Z}$, where the equality is easy to check.) A complex vector bundle has Chern classes in complex cobordism which map to the usual Chern classes in ordinary cohomology [2], so this equality means that $2\chi$ is in the image of $MU^4 BSO(4)$. For convenience, let $C$ be the element $c_2 A - c_2 B$ in $MU^4 BSO(4)$, so that $C$ maps to $2\chi$ in $H^4(BSO(4), \mathbf{Z})$.

To prove that $\chi$ itself is not in the image of $MU^4 BSO(4)$, we show that a certain odd-dimensional Steenrod operation, $Sq^3$, is nonzero on the image of $\chi$ in $H^4(BSO(4), \mathbf{Z}/2)$, which is a polynomial ring $\mathbf{Z}/2[w_2, w_3, w_4]$. We use Wu's formula for the Steenrod operations in $H^*(BO(n), \mathbf{Z}/2)$ [25], p. 94:

$$Sq^3 \chi = Sq^3 w_4 = w_3 w_4 \neq 0.$$

Thus $\chi \in H^4(BSO(4), \mathbf{Z})$ is not in the image of $MU^4 BSO(4)$. Since $H^4(BSO(4), \mathbf{Z}) = \mathbf{Z} \oplus \mathbf{Z}$ has no torsion, Corollary 5.3 applies, proving Theorem 5.1 for the case of $SO(4)$. (By the way, the same thing happens for $SO(2n)$ for all $n \geq 2$: $2^{n-1}\chi \in H^{2n}(BSO(2n), \mathbf{Z})$ is a polynomial in Chern classes of representations of $SO(2n)$, so it is in the image of $MU^{2n} BSO(2n)$, but $\chi$ itself is not in the image. It is plausible that $2^{n-1}\chi$ should be the smallest multiple of $\chi$ which is in the image of $MU^{2n} BSO(2n)$, as the above calculation shows for $SO(4)$ and as Inoue [17] showed for $SO(6)$.)

Now we turn to the construction of a similar example among finite groups. The idea is to use a finite subgroup $G \subset SO(4)$ and the restriction of the Euler class $\chi$ to $H^4(BG, \mathbf{Z})$. Then it is automatic that $2\chi$ is in the image of $MU^4 BG$, and we just have to choose $G$ so that $\chi + $ (2-torsion in $H^4(BG, \mathbf{Z})$) does not intersect the image of $MU^4 BG$. (Throughout this paper, a 2-torsion element of an abelian group will mean an element $x$ with $2x = 0$, not just an element killed by some power of 2.) Since the phenomenon we are considering is 2-local, it is natural to take $G$ to be a reasonably big 2-subgroup of $SO(4)$. In general, $SO(n)$ contains a fairly big abelian 2-subgroup, the group $(\mathbf{Z}/2)^{n-1}$ of diagonal matrices with entries $\pm 1$, but it turns out that abelian subgroups of $SO(4)$ do not have the property we want. Fortunately, we get a more interesting subgroup of $SO(4)$ by defining $G$ to be the inverse image of the subgroup $(\mathbf{Z}/2)^2 \times (\mathbf{Z}/2)^2 \subset SO(3) \times SO(3)$ under the double cover $SO(4) \to SO(3) \times SO(3)$. Thus $G$ is an extra-special group of order 32, that is, a central extension

$$1 \to \mathbf{Z}/2 \to G \to (\mathbf{Z}/2)^4 \to 1$$

with center exactly $\mathbf{Z}/2$. In different terminology, $G$ is a Heisenberg group and the embedding $G \subset SO(4)$, the unique irreducible representation of $G$ of dimension greater than one, is the Schrödinger representation of $G$.

In what is probably the most beautiful calculation in the cohomology of groups, Quillen [32] computed the $\mathbf{Z}/2$-cohomology of the extra-special 2-groups. For the above group $G$, which Quillen calls "the real case," where the central extension is classified by the quadratic form $x_1x_2+x_3x_4$ over $\mathbf{Z}/2$, Corollary 5.12 in Quillen's paper says that the embedding $G \subset SO(4)$ makes $H^*(BG, \mathbf{Z}/2)$ a free module over $H^*(BSO(4), \mathbf{Z}/2) = \mathbf{Z}/2[w_2, w_3, w_4]$. It follows that

$$Sq^3\chi = Sq^3 w_4 = w_3 w_4 \neq 0 \in H^*(BG, \mathbf{Z}/2),$$

since the same is true in $H^*(BSO(4), \mathbf{Z}/2)$. We need to show a little more than this, namely that if $y \in H^4(BG, \mathbf{Z})$ is killed by 2, then $Sq^3(\chi+y) \neq 0$. To see this, we have to use a second description of $H^*(BG, \mathbf{Z}/2)$ from Quillen's paper: $H^*(BG, \mathbf{Z}/2)$ is the tensor product of a quotient ring of $H^*(B(\mathbf{Z}/2)^4, \mathbf{Z}/2)$ with the polynomial ring $\mathbf{Z}/2[w_4]$. (The homomorphism $H^*(B(\mathbf{Z}/2)^4, \mathbf{Z}/2) \to H^*(BG, \mathbf{Z}/2)$ comes from the abelianization map $G \to (\mathbf{Z}/2)^4$.) Now if $y$ is an element of $H^4(BG, \mathbf{Z})$ killed by 2, then $y$ is the Bockstein of an element of $H^3(BG, \mathbf{Z}/2)$, and by Quillen's second description of $H^*(BG, \mathbf{Z}/2)$, all of $H^3(BG, \mathbf{Z}/2)$ is in the image of $H^*(B(\mathbf{Z}/2)^4, \mathbf{Z}/2)$. It follows that $Sq^3 y$ is also in the image of $H^*(B(\mathbf{Z}/2)^4, \mathbf{Z}/2)$. Since $Sq^3\chi = w_3 w_4$ is not in that subring, we have that $Sq^3(\chi+y) \neq 0$ for all $y \in H^4(BG, \mathbf{Z})$ killed by 2. This is what we need for Corollary 5.3 to apply. Thus the proof of Theorem 5.1 is complete. In particular, the map $MU^6(BG \times B\mathbf{Z}/2) \otimes_{MU^*} \mathbf{Z} \to H^6(BG \times B\mathbf{Z}/2, \mathbf{Z})$ is not injective. QED

**6. Finite complexes with $MU^*X \otimes_{MU^*} \mathbf{Z} \to H^*(X, \mathbf{Z})$ not injective.** In section 5, we gave examples of compact Lie groups $G$ such that the maps $MU^*BG \otimes_{MU^*} \mathbf{Z}/2 \to H^*(BG, \mathbf{Z}/2)$ or $MU^*BG \otimes_{MU^*} \mathbf{Z} \to H^*(BG, \mathbf{Z})$ are not injective. In this section, we show that the elements we construct in the kernel remain nonzero when restricted from $BG$ to its $n$-skeleton, for some finite dimension $n$. This is not hard to prove if we don't need to know exactly what dimension is necessary, but we prefer to prove this for the smallest possible dimension, although this seems to require a long calculation. The result will be used in section 7 to construct our examples in algebraic geometry.

Let $G$ be the central extension
$$1 \to \mathbf{Z}/2 \to G \to (\mathbf{Z}/2)^4 \to 1$$

considered in section 5. We defined an element $C \in MU^4 BG$ which is nonzero in $MU^4 BG \otimes_{MU^*} \mathbf{Z}/2$ but which maps to 0 in $H^4(BG, \mathbf{Z}/2)$.

**Proposition 6.1** *The element $C \in MU^4 BG \otimes_{MU^*} \mathbf{Z}/2$ remains nonzero in $MU^4 X_7 \otimes_{MU^*} \mathbf{Z}/2$, where $X_7$ denotes the 7-skeleton of (any cell decomposition of) $BG$.*

**Proof.** We proved that $C \in MU^4 BG$ is nonzero in $MU^4 BG \otimes_{MU^*} \mathbf{Z}/2$ by computing that $C$ maps to $2\chi \in H^4(BG, \mathbf{Z})$, where $Sq^3(\chi +$ (any 2-torsion element in $H^4(BG, \mathbf{Z})))$ is nonzero in $H^7(BG, \mathbf{Z})$. The same calculations apply to the 7-skeleton of $BG$. QED

We now consider the element $C \otimes c_1 \in MU^6(BG \times B\mathbf{Z}/2) \otimes_{MU^*} \mathbf{Z}$, for $G$ the central extension mentioned above, which maps to 0 in $H^6(BG \times B\mathbf{Z}/2, \mathbf{Z})$ by section 5.

**Proposition 6.2** *The element $C \otimes c_1 \in MU^6(BG \times B\mathbf{Z}/2) \otimes_{MU^*} \mathbf{Z}$ remains nonzero in $MU^6(BG \times B\mathbf{Z}/2)_{15} \otimes_{MU^*} \mathbf{Z}$, where $(BG \times B\mathbf{Z}/2)_{15}$ denotes the 15-skeleton of $BG \times B\mathbf{Z}/2$.*

The proof occupies the rest of this section.
**Proof.** Slightly more precisely, we will show that $C \otimes c_1$ is nonzero in $MU^6(X_7 \times Y_8) \otimes_{MU^*} \mathbf{Z}$, where $X_7$ is the 7-skeleton of $BG$ and $Y_8$ is the 8-skeleton of $B\mathbf{Z}/2$.

The space $B\mathbf{Z}/2$ is the $S^1$-bundle over $\mathbf{CP}^\infty$ associated to the complex line bundle $L^{\otimes 2}$, where $L$ is the hyperplane line bundle on $\mathbf{CP}^\infty$. We have $BP^*\mathbf{CP}^\infty \cong BP^*[[c_1]]$. As Stong [38] first observed, the Gysin sequence for this $S^1$-bundle gives that $BP^*B\mathbf{Z}/2 = BP^*[[c_1]]/([2](c_1) = 0)$, where $[2](c_1)$ is the power series with coefficients in the ring $BP^*$ which computes $c_1(L^{\otimes 2})$ in terms of $c_1 = c_1(L) \in BP^2\mathbf{CP}^\infty$. We adopt the convention that $v_i \in BP^{-2(2^i-1)}$ means the coefficient of $c_1^{2^i}$ in this power series; these coefficients are polynomial generators for the ring $BP^*$ [**42**], p. 20. We will only need to know the first few terms of the series $[2](c_1)$ in terms of these generators:

$$[2](c_1) = 2c_1 + v_1 c_1^2 + 2v_1^2 c_1^3 + v_2 c_1^4 + \cdots$$

[**42**], p. 21.

Since the 8-skeleton $Y_8 \subset B\mathbf{Z}/2$ has $BP^i Y_8 = 0$ for $i > 8$, we have $c_1^5 = 0 \in BP^* Y_8$. Moreover, the cohomology ring $H^*(Y_8, \mathbf{Z}) = \mathbf{Z}[c_1]/(2c_1 = 0, c_1^5 = 0)$ is generated by $c_1$, so the map $BP^* Y_8 \to H^*(Y_8, \mathbf{Z}_{(2)})$ is surjective. It follows that the Atiyah-Hirzebruch spectral sequence $H^*(Y_8, BP^*) \Rightarrow BP^* Y_8$ degenerates, and we can then read off that the natural map

$$BP^*[[c_1]]/(2c_1 + v_1 c_1^2 + 2v_1^2 c_1^3 + v_2 c_1^4 = 0, c_1^5 = 0) \to BP^* Y_8$$

is an isomorphism. The same argument determines $BP^* Y_n$ for all even $n$, and we will need that $BP^* Y_2 = BP^*[[c_1]]/(2c_1 = 0, c_1^2 = 0)$.

We can now begin to describe $BP^*(X_7 \times Y_8)$.

**Lemma 6.3** *Suppose $X$, $Y$ are finite cell complexes such that $BP^* Y \to H^*(Y, \mathbf{Z}_{(p)})$ is surjective, where $BP$ denotes the Brown-Peterson cohomology theory associated to a prime number $p$. Then there is an exact sequence of $BP^*$-modules,*

$$0 \to BP^* X \otimes_{BP^*} BP^* Y \to BP^*(X \times Y) \to \operatorname{Tor}_1^{BP^*}(BP^* X, BP^* Y) \to 0.$$

**Proof.** By Johnson-Wilson [**18**], for a finite complex $Y$, surjectivity of $BP^* Y \to H^*(Y, \mathbf{Z}_{(p)})$ is equivalent to $BP^* Y$'s having projective dimension $\leq 1$ as a $BP^*$-module. In particular, $\operatorname{Tor}_i^{BP^*}(BP^* X, BP^* Y) = 0$ for $i \geq 2$. Given that, the Künneth spectral sequence for $BP$-theory [**1**] reduces to the above exact sequence. QED

We apply this to $X_7 \subset BG$ and $Y_{2n} \subset B\mathbf{Z}/2$. We have already mentioned that $BP^* Y_{2n} \to H^*(Y_{2n}, \mathbf{Z}_{(2)})$ is surjective, so we have a short exact sequence of $BP^*$-modules,

$$0 \to BP^* X_7 \otimes_{BP^*} BP^* Y_{2n} \to BP^*(X_7 \times Y_{2n}) \to \operatorname{Tor}_1^{BP^*}(BP^* X_7, BP^* Y_{2n}) \to 0.$$

We want to show that the element $C \otimes c_1 \in BP^4 X_7 \otimes_{\mathbf{Z}} BP^2 Y_8 \to BP^6(X_7 \times Y_8)$ is nonzero in $BP^6(X_7 \times Y_8) \otimes_{BP^*} \mathbf{Z}$, that is, that we cannot write

$$C \otimes c_1 = \sum v_i x_i$$

with $x_i \in BP^{6+2(2^i-1)}(X_7 \times Y_8)$. Thus, $x_1 \in BP^8(X_7 \times Y_8)$, $x_2 \in BP^{12}(X_7 \times Y_8)$, $x_3 \in BP^{20}(X_7 \times Y_8)$, and so on. (Since $X_7 \times Y_8$ only has dimension 15, the elements $x_i$ for $i \geq 3$ are 0.) The idea here is that in the limit $Y_\infty = B\mathbf{Z}/2$, the Tor term in the above exact sequence is 0 (that is, $BP^* B\mathbf{Z}/2$ is flat in the sense of the proof of Lemma 5.2). So if we consider sufficiently high-dimensional skeleta of $B\mathbf{Z}/2$, the 8-skeleton being enough for our purpose, any complications caused by that Tor term will go away.

We will prove that we cannot write $C \otimes c_1 = \sum v_i x_i$ in $BP^6(X_7 \times Y_8)$ by comparing $X_7 \times Y_8$ with its subspace $X_7 \times Y_2$. Consider the map of exact sequences

$$\begin{array}{ccccccccc}
0 & \to & BP^* X_7 \otimes_{BP^*} BP^* Y_8 & \to & BP^*(X_7 \times Y_8) & \to & \operatorname{Tor}_1^{BP^*}(BP^* X_7, BP^* Y_8) & \to & 0 \\
& & \downarrow & & \downarrow & & \downarrow & & \\
0 & \to & BP^* X_7 \otimes_{BP^*} BP^* Y_2 & \to & BP^*(X_7 \times Y_2) & \to & \operatorname{Tor}_1^{BP^*}(BP^* X_7, BP^* Y_2) & \to & 0.
\end{array}$$

If we have the equation $C \otimes c_1 = \sum v_i x_i$ in $BP^6(X_7 \times Y_8)$, as above, then we can restrict it to get an equality in $BP^6(X_7 \times Y_2)$. Since $X_7 \times Y_2$ has dimension only 9, the elements $x_i$ restrict to 0 for $i \geq 2$. Suppose we can show that the image of $x_1$ in $\operatorname{Tor}_1^{BP^*}(BP^* X_7, BP^* Y_8)$ maps to 0 in $\operatorname{Tor}_1^{BP^*}(BP^* X_7, BP^* Y_2)$. Then we have $C \otimes c_1 = v_1 x_1$ for some $x_1 \in BP^* X_7 \otimes_{BP^*} BP^* Y_2$, so that $C \otimes c_1$ is 0 in $(BP^* X_7 \otimes_{BP^*} BP^* Y_2) \otimes_{BP^*} \mathbf{Z}_{(2)}$.

But we know (Proposition 6.1) that $C$ is nonzero in $BP^*X_7 \otimes_{BP^*} \mathbf{Z}/2$, and $c_1$ generates a $\mathbf{Z}/2$ summand of the abelian group $BP^*Y_2 \otimes_{BP^*} \mathbf{Z}_{(2)} = H^*(Y_2, \mathbf{Z}_{(2)}) = \mathbf{Z}_{(2)} \oplus \mathbf{Z}/2$, so that $C \otimes c_1$ is nonzero in $BP^*X_7 \otimes_{BP^*} (BP^*Y_2 \otimes_{BP^*} \mathbf{Z}_{(2)}) = (BP^*X_7 \otimes_{BP^*} BP^*Y_2) \otimes_{BP^*} \mathbf{Z}_{(2)}$. This contradiction shows that, if we can show that $x_1$ maps to 0 in $\text{Tor}_1^{BP^*}(BP^*X_7, BP^*Y_2)$, then $C \otimes c_1$ is nonzero in $BP^6(X_7 \times Y_8) \otimes_{BP^*} \mathbf{Z}_{(2)}$, proving Proposition 6.2.

The fact we need is supplied by the following lemma.

**Lemma 6.4** *Any element $x_1 \in \text{Tor}_1^{BP^*}(BP^*X_7, BP^*Y_8)$ of degree 8 such that $\sum_{i \geq 1} v_i x_i = 0 \in \text{Tor}_1^{BP^*}(BP^*X_7, BP^*Y_8)$ for some $x_2, x_3, \ldots$ must restrict to 0 in $\text{Tor}_1^{BP^*}(BP^*X_7, BP^*Y_2)$.*

This applies to the class $x_1$ considered above. In fact, we assumed an equality $C \otimes c_1 = \sum v_i x_i$ in $BP^6(X_7 \times Y_8)$, where $C \otimes c_1$ lies in $BP^*X_7 \otimes_{BP^*} BP^*Y_8 \subset BP^*(X_7 \times Y_8)$. So when you map this equation into $\text{Tor}_1^{BP^*}(BP^*X_7, BP^*Y_8)$, it becomes $0 = \sum v_i x_i$.

**Proof of Lemma 6.4.** We begin by describing $\text{Tor}_1^{BP^*}(BP^*X_7, BP^*Y_{2n})$ for any $n$. We can give an explicit free resolution of $BP^*Y_{2n} = BP^*[[c_1]]/(2c_1 + v_1c_1^2 + 2v_1^2c_1^3 + v_2c_1^4 + \cdots = 0, c_1^{2(n+1)} = 0)$ as a $BP^*$-module:
$$0 \to \oplus_{\text{degrees } 2,4,\ldots,2n} BP^* \to \oplus_{\text{degrees } 0,2,4,\ldots,2n} BP^* \to BP^*Y_{2n} \to 0.$$
That is, $BP^*Y_{2n}$ is generated as a $BP^*$-module by $1, c_1, c_1^2, \ldots, c_1^n$, with relations $2c_1^n = 0, 2c_1^{n-1} + v_1c_1^n = 0, \ldots, 2c_1 + v_1c_1^2 + 2v_1^2c_1^3 + v_2c_1^4 + \cdots = 0$. So
$$\text{Tor}_1^{BP^*}(BP^*X_7, BP^*Y_{2n}) \cong \{(x_1, x_2, \ldots, x_n) \in (BP^*X_7)^n : 2x_1 = 0, 2x_2 + v_1x_1 = 0,$$
$$\ldots, 2x_n + v_1x_{n-1} + 2v_1^2 x_{n-2} + v_2 x_{n-3} + \cdots = 0\},$$
graded so that elements of degree $i$ have $x_1 \in BP^{i-1}X_7$, $x_2 \in BP^{i-3}X_7$, $x_3 \in BP^{i-5}X_7$, and so on. Finally, the map $\text{Tor}_1^{BP^*}(BP^*X_7, BP^*Y_8) \to \text{Tor}_1^{BP^*}(BP^*X_7, BP^*Y_2)$ sends a quadruple $x_1, x_2, x_3, x_4$ of elements of $BP^*X_7$ satisfying these equations to the first element $x_1$.

In the specific situation of the lemma, the group $\text{Tor}_1^{BP^*}(BP^*X_7, BP^*Y_8)$ in degree 8 is the group of quadruples $x_1 \in BP^7X_7$, $x_2 \in BP^5X_7$, $x_3 \in BP^3X_7$, $x_4 \in BP^1X_7$ which satisfy the equations $2x_1 = 0$, $2x_2 + v_1x_1 = 0$, $2x_3 + v_1x_2 + 2v_1^2x_1 = 0$, and $2x_4 + v_1x_3 + 2v_1^2x_2 + v_2x_1 = 0$. We have to show that if we can write $\sum_{i \geq 1} v_i(x_1^i, \ldots, x_4^i) = 0$ for some quadruples $(x_1^i, \ldots, x_4^i)$ satisfying these equations, with $(x_1^1, \ldots, x_4^1)$ the given elements $(x_1, \ldots, x_4)$, then $x_1 = 0 \in BP^7X_7$.

For this we need to know something about the structure of $BP^*X_7$ as a $BP^*$-module. Most of what we need follows just from $X_7$ having low dimension, namely 7. Consider the Atiyah-Hirzebruch spectral sequence $H^*(X_7, BP^*) \Longrightarrow BP^*X_7$, which is in the lower right quadrant. We write $H^iX_7$ for $H^i(X_7, \mathbf{Z}_{(2)})$ in this spectral sequence.

| $H^0X_7$ | $H^1X_7$ | $H^2X_7$ | $H^3X_7$ | $H^4X_7$ | $H^5X_7$ | $H^6X_7$ | $H^7X_7$ |
|---|---|---|---|---|---|---|---|
| $H^0X_7 \cdot v_1$ | $H^1X_7 \cdot v_1$ | $H^2X_7 \cdot v_1$ | $H^3X_7 \cdot v_1$ | $H^4X_7 \cdot v_1$ | $H^5X_7 \cdot v_1$ | $H^6X_7 \cdot v_1$ | $H^7X_7 \cdot v_1$ |
| $H^0X_7 \cdot v_1^2$ | $H^1X_7 \cdot v_1^2$ | $H^2X_7 \cdot v_1^2$ | $H^3X_7 \cdot v_1^2$ | $H^4X_7 \cdot v_1^2$ | $H^5X_7 \cdot v_1^2$ | $H^6X_7 \cdot v_1^2$ | $H^7X_7 \cdot v_1^2$ |
| $H^0X_7 \cdot v_1^3$ | $H^1X_7 \cdot v_1^3$ | $H^2X_7 \cdot v_1^3$ | $H^3X_7 \cdot v_1^3$ | $H^4X_7 \cdot v_1^3$ | $H^5X_7 \cdot v_1^3$ | $H^6X_7 \cdot v_1^3$ | $H^7X_7 \cdot v_1^3$ |
| $\oplus$ | $\oplus$ | $\oplus$ | $\oplus$ | $\oplus$ | $\oplus$ | $\oplus$ | $\oplus$ |
| $H^0X_7 \cdot v_2$ | $H^1X_7 \cdot v_2$ | $H^2X_7 \cdot v_2$ | $H^3X_7 \cdot v_2$ | $H^4X_7 \cdot v_2$ | $H^5X_7 \cdot v_2$ | $H^6X_7 \cdot v_2$ | $H^7X_7 \cdot v_2$ |

$\vdots$

The odd rows are 0 and so are not written in the above diagram. The $E_2$ term is generated as an algebra by $H^*(X_7, \mathbf{Z}_{(2)})$ in the top row and $BP^*$ in the left column. All the differentials are 0 on $BP^*$ in the left column (because $BP^*X_7$ maps onto $BP^*(\text{base point})$), so it suffices to describe the first nonzero differential, $d_3$, on $H^*(X_7, \mathbf{Z}_{(2)})$ in the top row. This first differential, from $H^*(X_7, \mathbf{Z}_{(2)})$ to $H^*(X_7, \mathbf{Z}_{(2)}) \cdot v_1$, is known to be the Steenrod operation $Sq^3$, which makes sense on integral cohomology by defining it as the composition

$$H^i(X, \mathbf{Z}) \to H^i(X, \mathbf{Z}/2) \xrightarrow[Sq^2]{} H^{i+2}(X, \mathbf{Z}/2) \xrightarrow[\beta]{} H^{i+3}(X, \mathbf{Z}).$$

After that first differential, what is left in the spectral sequence is

$$(\ker Sq^3 \subset H^*(X_7, \mathbf{Z}_{(2)})) \otimes_{\mathbf{Z}_{(2)}} BP^* / ((\operatorname{im} Sq^3 \subset H^*(X_7, \mathbf{Z}_{(2)})) \cdot v_1 = 0).$$

In particular, the spectral sequence is still generated as an algebra by the top row, $\ker Sq^3 \subset H^*(X_7, \mathbf{Z}_{(2)})$, and $BP^*$ in the left column. But then, since $\dim X_7 \leq 7$, there are no more differentials: any other differentials would have horizontal degree at least 5, so they would be defined on $H^i(X_7, \mathbf{Z}_{(2)})$ for $i \leq 2$, but we know (Theorem 2.2) that $BP^*X \to H^*(X, \mathbf{Z}_{(2)})$ is always surjective for $i \leq 2$ so that there are no nonzero differentials on $H^i(X, \mathbf{Z}_{(2)})$ for $i \leq 2$. Thus the $E_\infty$ term of the spectral sequence converging to $BP^*X_7$ is $(\ker Sq^3 \subset H^*(X_7, \mathbf{Z}_{(2)})) \otimes_{\mathbf{Z}_{(2)}} BP^* / ((\operatorname{im} Sq^3 \subset H^*(X_7, \mathbf{Z}_{(2)})) \cdot v_1 = 0)$. In particular, multiplication by $v_i$, $i \geq 2$, is injective on the $E_\infty$ term. Also, multiplication by $v_1$ is injective except on $H^6(X_7, \mathbf{Z}_{(2)})$ or $H^7(X_7, \mathbf{Z}_{(2)})$ times monomials in $v_2, v_3, \ldots$, because $Sq^3$ is at most nonzero from $H^3$ to $H^6$ or from $H^4$ to $H^7$.

This gives a lot of information on the elements $x_1 = x_1^1 \in BP^7 X_7, \ldots, x_4 = x_4^1 \in BP^1 X_7$, just using that they satisfy equations $\sum v_i x_r^i = 0$, $r = 1, \ldots, 4$. In particular, consider what this says about $x_4 \in BP^1 X_7$. Each group $BP^i X_7$ is filtered with quotients certain groups in the $E_\infty$ term in the Atiyah-Hirzebruch spectral sequence, so any nonzero element of $BP^i X_7$ has a "leading term" which is a nonzero element of one of the groups in the $E_\infty$ term. The leading term of $x_4$ cannot be in $H^1(X_7, \mathbf{Z}_{(2)})$, since multiplication by $v_1$ is an injection from that group into row $-2$ of the spectral sequence, whereas the sum $\sum_{i \geq 2} v_i x_4^i$, if nonzero, will have leading term in some row $\leq -6$; in fact, if its leading term is anywhere above row $-14$ (where $v_3$ is), it must be in some group times $v_2$. Likewise the leading term of $x_4$ cannot be in $(\ker Sq^3 \subset H^3(X_7, \mathbf{Z}_{(2)})) \cdot v_1$ or in $H^5(X_7, \mathbf{Z}_{(2)}) \cdot v_1^2$ or even in the $H^7(X_7, \mathbf{Z}_{(2)})/Sq^3 \cdot v_1^3$ summand in row $-6$: either $x_4 \in BP^1 X_7$ is 0 or its leading term is in $H^7(X_7, \mathbf{Z}_{(2)}) \cdot v_2 \subset BP^1 X_7$. The same argument shows that $x_3 \in BP^3 X_7$ and $x_2 \in BP^5 X_7$ must actually be 0.

Now we use the other equations satisfied by $x_1, \ldots, x_4$: $2x_1 = 0$, $2x_2 + v_1 x_1 = 0$, $2x_3 + v_1 x_2 + 2v_1^2 x_1 = 0$, and $2x_4 + v_1 x_3 + 2v_1^2 x_2 + v_2 x_1 = 0$. Since $x_2 = x_3 = 0$, the last equation becomes $2x_4 + v_2 x_1 = 0$. Here $x_1 \in BP^7 X_7$, and there is only one group in the spectral sequence which contributes to $BP^7 X_7$, namely $BP^7 X_7 = H^7(X_7, \mathbf{Z}_{(2)})$ on the top row. Since multiplication by $v_2$ is injective on the spectral sequence, the equation $2x_4 + v_2 x_1 = 0$, with $x_4 \in H^7(X_7, \mathbf{Z}_{(2)}) \cdot v_2$, implies that $x_1 \in H^7(X_7, \mathbf{Z}_{(2)})$ is a multiple of 2.

Now it is convenient to use Harada and Kono's calculation of the integral cohomology of $BG$ ($X_7$ above being the 7-skeleton of $BG$), going beyond Quillen's computation of the $\mathbf{Z}/2$-cohomology which I used earlier [**14**]. They showed in particular that $H^i(BG, \mathbf{Z})$ is a $\mathbf{Z}/2$-vector space for $i \not\equiv 0 \pmod 4$. Now the 7-skeleton $X_7$ has $H^7(X_7, \mathbf{Z}) = H^7(BG, \mathbf{Z}) \oplus$ (free abelian group), so $H^7(X_7, \mathbf{Z})$ also has no 4-torsion. Since $2x_1 = 0$ and $x_1$ is a multiple of 2, we have $x_1 = 0 \in H^7(X_7, \mathbf{Z}_{(2)}) \cong BP^7 X_7$. This is what we needed to prove Lemma 6.4. QED

## 7. Non-injectivity of the classical cycle maps in algebraic geometry.

**Theorem 7.1** *There is a smooth complex projective variety $X$ of dimension 7 such that the map $CH^2 X/2 \to H^4(X, \mathbf{Z}/2)$ is not injective. The variety $X$ and the element of $CH^2 X/2$ in the kernel which we construct can be defined over $\mathbf{Q}$.*

**Proof.** Let $G$ be the Heisenberg group

$$1 \to \mathbf{Z}/2 \to G \to (\mathbf{Z}/2)^4 \to 1,$$

as in section 5. Let $X = Y/G$ be the quotient of a complete intersection $Y \subset P(V)$ by a free action of $G$, where $G$ acts linearly on the vector space $V$. Such varieties exist for any finite group $G$, over any infinite field, and for $X$ of any dimension $\geq 1$, by Godeaux and Serre; see [**36**], section 20. Since $X$ is the quotient of a free $G$-action, there is a natural homotopy class of maps $X \to BG$, and since the $G$-action on $Y$ is "linearized" there is a natural homotopy class of maps $X \to \mathbf{CP}^\infty$ (or equivalently, a natural generator of $H^2(X, \mathbf{Z})$). By Atiyah and Hirzebruch [4], p. 42, the product map $X \to BG \times \mathbf{CP}^\infty$ is $r$-connected. In particular, $X$ contains the $r$-skeleton of $BG$ up to homotopy.

In section 5, we defined two representations $A$ and $B$ of $SO(4)$, of dimensions 3 and 4. We restrict these to $G \subset SO(4)$, and we define $C = c_2 A - c_2 B \in MU^4 BG$. We proved that $C$ maps to 0 in $H^4(BG, \mathbf{Z}/2)$, but $C$ remains nonzero in $MU^4 X_7 \otimes_{MU^*} \mathbf{Z}/2$, where $X_7$ denotes the 7-skeleton of $BG$, by Proposition 6.1. Let $X = Y/G$ be a Godeaux-Serre variety for this group with complex dimension at least 7. Since $X$ contains the 7-skeleton of $BG$ up to homotopy, the class $C \in MU^4 BG$ pulls back to a nonzero element of $MU^4 X \otimes_{MU^*} \mathbf{Z}/2$, and it clearly maps to 0 in $H^4(X, \mathbf{Z}/2)$.

Moreover, the complex representations $A$ and $B$ of $G$ give algebraic vector bundles over $X = Y/G$, and we can consider the algebraic cycle $C := c_2 A - c_2 B \in CH^2 X$. It maps to the above class $C \in MU^4 X \otimes_{MU^*} \mathbf{Z}/2$, which implies that $C$ is nonzero in $CH^2 X/2$ but maps to 0 in $H^4(X, \mathbf{Z}/2)$.

The variety $X$ can be defined over $\mathbf{Q}$ (or any infinite field) by Serre's construction, and the cycle $C \in CH^2 X$ can be defined over $\mathbf{Q}$ because the representations $A$ and $B$ of the group $G$ can be defined over $\mathbf{Q}$. QED

**Theorem 7.2** *There is a smooth complex projective variety $X$ of dimension 15 and an element $\alpha \in CH^3 X$ with the following properties:*
$2\alpha = 0 \in CH^3 X$;
$\alpha$ *maps to 0 in $H^6(X, \mathbf{Z})$ and also in the intermediate Jacobian $H^5(X, \mathbf{C})/(F^3 H^5(X, \mathbf{C}) + H^5(X, \mathbf{Z}))$;*
$\alpha$ *is not algebraically equivalent to 0.*

**Proof.** Let $G$ be the Heisenberg group

$$1 \to \mathbf{Z}/2 \to G \to (\mathbf{Z}/2)^4 \to 1$$

as above, and let $X$ be a Godeaux-Serre variety $X = Y/(G \times \mathbf{Z}/2)$ of dimension at least 15. Let $\alpha$ be the class $Cc_1 \in CH^3 X$. Here $C = c_2 A - c_2 B$, where $A$ and $B$ are the 3- and 4-dimensional representations of $G$ considered above, and $c_1$ denotes the first Chern class of the nontrivial character of $\mathbf{Z}/2$. Clearly $2\alpha = 0$, since $2c_1 = 0$.

The image of $\alpha$ in $MU^6 X \otimes_{MU^*} \mathbf{Z}$ is the pullback to $X$ of the element $C \otimes c_1 \in MU^6(BG \times B\mathbf{Z}/2) \otimes_{MU^*} \mathbf{Z}$ considered in Proposition 6.2. Since $X$ contains the 15-skeleton of $BG \times B\mathbf{Z}/2$ up to homotopy, that lemma implies that $\alpha$ is nonzero in $MU^6 X \otimes_{MU^*} \mathbf{Z}$ but maps to 0 in $H^6(X, \mathbf{Z})$.

Since $\alpha$ is nonzero in $MU^6 X \otimes_{MU^*} \mathbf{Z}$, it is not algebraically equivalent to 0. The intermediate Jacobian for codimension-3 cycles on $X$, $H^5(X, \mathbf{C})/(F^3 H^5(X, \mathbf{C}) + H^5(X, \mathbf{Z}))$, is actually 0, since $H^5(X, \mathbf{C}) \subset H^5(Y, \mathbf{C})$ and $Y$ is a complete intersection of dimension $\geq 15$. Thus $\alpha$ is 0 in the intermediate Jacobian as well as in ordinary cohomology. QED

**8. Further comments. Remark 1.** Our main example is a codimension-3 cycle on a smooth projective variety of rather large dimension, 15. It is worth mentioning that in a sense this large dimension should be inessential. Namely, if $X$ is a smooth complex projective variety and $Y$ is a smooth ample hypersurface in $X$, then the restriction map $CH^i X \to CH^i Y$ is conjectured to be an isomorphism for $i < \dim Y/2$, which would be a version of the Lefschetz hyperplane theorem for Chow groups [15], [28], p. 643. Moreover if $Y$ is a very general smooth hypersurface whose class in $CH^1 X$ is a sufficiently high multiple of an ample class, then Nori conjectured that much more should be true: $CH^i X \to CH^i Y$ should be an isomorphism for all $i < \dim Y$ [26], p. 368, [28], p. 644. Actually Nori and Paranjape only state these conjectures $\otimes \mathbf{Q}$, but they seem plausible integrally in view of Kollár and van Geemen's Trento examples [5], p. 135. (Nori also conjectured that for $Y$ a very general high-degree complete intersection in a smooth projective variety $X$, $CH^i X \otimes \mathbf{Q} \to CH^i Y \otimes \mathbf{Q}$ should be injective for $i = \dim Y$. Here we cannot expect to have the corresponding integral statement: if $\alpha \in CH^3 X$ is the cycle in Theorem 7.2, then $\alpha$ restricts to 0 in $CH^3 Y$ for every complete intersection 3-fold $Y \subset X$, by Roitman's theorem.)

These conjectures would imply the corresponding isomorphisms for cycles modulo algebraic equivalence in place of Chow groups. In particular, the integral version of Nori's conjecture would imply that our nonzero element of the Griffiths group ker $(Z^3_{\text{alg}} X \to H^6(X, \mathbf{Z}))$, for dim $X = 15$, remains nonzero on a very general high-degree complete intersection $X' \subset X$ of dimension as small as 4.

But we could not expect to prove that our cycle remains nonzero in $Z^3_{\text{alg}} X'$ for such a small-dimensional variety $X'$ just using the cycle class defined in this paper. In fact, on a variety $X'$ of such small dimension, our cycle would be 0 in $MU^6 X' \otimes_{MU^*} \mathbf{Z}$ as well as in $H^6(X', \mathbf{Z})$, because for any finite cell complex $X$ of real dimension at most $14 = 2(2^2 + 2 + 1)$, the map $MU^* X \otimes_{MU^*} \mathbf{Z} \to H^*(X, \mathbf{Z})$ is injective. (Proof: By Johnson and Wilson [18], Proposition 4.1 and Theorem 1.1, dim $X \leq 14$ implies that hom dim$_{MU^*} MU^* X \leq 2$, and by Conner and Smith [11], that implies that $MU^* X \otimes_{MU^*} \mathbf{Z} \to H^*(X, \mathbf{Z})$ is injective.)

**Remark 2.** Bloch defined an interesting filtration of the group of algebraic cycles homologically equivalent to 0, with the smallest subgroup being the cycles algebraically equivalent to 0 [7], p. 380. Namely, one says that a $k$-dimensional cycle $\alpha$ on a variety $X$ is $r$-equivalent to 0 if $\alpha$ is contained in some $(k+r)$-dimensional algebraic subset $S \subset X$ such that $[\alpha] = 0 \in H^{\text{BM}}_{2k}(S, \mathbf{Z})$. Then 1-equivalence is the same thing as algebraic equivalence by Bloch [7], Lemma 1.1 (stated $\otimes \mathbf{Q}$, but the proof works integrally), and $r$-equivalence is the same as homological equivalence for $r \geq \dim X - \dim \alpha$. The cycle map $Z^{\text{alg}}_* X \to MU^{\text{BM}}_* X \otimes_{MU_*} \mathbf{Z}$ is well-defined on 1-equivalence by Theorem 3.1, and it is not well-defined on 3-equivalence by Theorem 7.2, since that result gives a codimension-3 cycle which is homologically equivalent to 0 but nonzero in $MU^{\text{BM}}_* X \otimes_{MU_*} \mathbf{Z}$. It seems possible that the cycle map is well-defined on 2-equivalence, but I can only prove a weaker statement, as follows.

Define a $k$-dimensional cycle $\alpha$ on a variety $X$ to be strongly $r$-equivalent to 0 if $\alpha$ is the pushforward, as a cycle, of a cycle $\alpha'$ on some *smooth* scheme $S'$ of dimension $\leq k + r$ with a proper map $S' \to X$, such that $[\alpha'] = 0 \in H^{\text{BM}}_{2k}(S', \mathbf{Z})$. Then strong 1-equivalence is the same as 1-equivalence, i.e., algebraic equivalence, by the proof of Lemma 1.1 in [7]; for $r \geq 2$ $r$-equivalence is in general different from strong $r$-equivalence, although they are the same $\otimes \mathbf{Q}$ under the assumption of the Hodge conjecture (and 2-equivalence $\otimes \mathbf{Q}$ is the same as strong 2-equivalence $\otimes \mathbf{Q}$ without any conjecture). Anyway, the point of this definition is that we can prove that the cycle map is well-defined on strong 2-equivalence, as follows. If $\alpha$ is the pushforward of a cycle $\alpha'$ which is homologically equivalent to 0 on a smooth $(k+2)$-dimensional scheme $S'$, then $[\alpha'] = 0$ in $MU^{\text{BM}}_{2k} S \otimes_{MU_*} \mathbf{Z}$ by the extension of Quillen's theorem given in Theorem 2.2, where we use injectivity in degree 4. It follows that $[\alpha] = 0 \in MU^{\text{BM}}_{2k} X \otimes_{MU_*} \mathbf{Z}$. That is, the cycle map is well-defined on strong 2-equivalence.

Department of Mathematics, University of Chicago, 5734 S. University Ave., Chicago, IL 60637
totaro@math.uchicago.edu